\documentclass[12pt]{iopart}
\usepackage{tikz}
\usepackage{pgf}
  \usepackage{color}
  \definecolor{purple}{rgb}{0.5,0,1}
  \definecolor{orange}{cmyk}{0,0.7,1,0}
  \definecolor{midgrey}{gray}{0.5}
\usepackage[english]{babel}
\usepackage{amsfonts}
\usepackage{amssymb}
\usepackage{amsthm}
\usepackage[numbers]{natbib}

\newcommand{\mbb}[1]{\mathbb{#1}}
\newcommand{\mc}[1]{\mathcal{#1}}
\newcommand{\binom}[2]{\left(\!\!\!\begin{array}{c}{#1}\\{#2}\end{array}\!\!\!\right)}

\eqnobysec
\newtheorem{theorem}{\textbf{Theorem}}[section]
\newtheorem{corollary}[theorem]{\textbf{Corollary}}
\newtheorem{lemma}[theorem]{Lemma}

\newtheorem{proposition}{\textbf{Proposition}}[section]

\newcommand{\Rea}{\mathrm{Re}\,}
\newcommand{\Ima}{\mathrm{Im}\,}
\newcommand{\sgn}{\mathrm{sgn}\,}
\newcommand{\sg}{\mathrm{sg}\,\alpha}

\begin{document}
\title[Solutions to some fragmentation equations with growth or decay]{Explicit solutions to some
fragmentation equations with growth or decay\footnote{The research has been partially supported by the National
Science Centre of Poland Grant 2017/25/B/ST1/00051 and the National Research Foundation of South Africa Grant 82770} }

\author{Jacek Banasiak$^{1,2}$,  David Wetsi Poka$^1$ and Sergey Shindin$^3$}
\address{$^1$ Department of Mathematics and Applied Mathematics, University of Pretoria, Pretoria, South Africa}
\address{$^2$ Institute of Mathematics, Technical University of \L\'{o}d\'z, \L\'{o}d\'z, Poland}
\address{$^3$ School of Mathematics, Statistics and Computer Science, University of KwaZulu-Natal, Durban, South Africa}
\eads{\mailto{jacek.banasiak@up.ac.za}, \mailto{wetsi@aims.ac.za},\mailto{shindins@ukzn.ac.za}}
\vspace{10pt}
\begin{indented}
\item[]We dedicate this paper to Robert M. Ziff in honour of his 70th birthday
\end{indented}

\begin{abstract}
In this paper, we provide a systematic way of finding explicit solutions for a class of continuous fragmentation equations with growth or decay in the state space and derive new explicit solutions in the cases of constant and linear growth/decay coefficients.
\end{abstract}

\section{Introduction}\label{sec1}
\begin{figure}
\begin{center}
\begin{tikzpicture}[scale=1]
\draw [cyan,fill]  (-2,0) circle [radius=0.15];
\draw [cyan,fill]  (-2,0.3) circle [radius=0.15];
\draw [cyan,fill]  (-2,-0.3) circle [radius=0.15];
\draw [red,fill]  (-1.4,0) circle [radius=0.15];
\draw [red,fill]  (-1.4,0.3) circle [radius=0.15];
\draw [red,fill]  (-1.4,-0.3) circle [radius=0.15];
\draw [yellow,fill]  (-1.7,0) circle [radius=0.15];
\draw [cyan,fill]  (2,0.8) circle [radius=0.15];
\draw [cyan,fill]  (2,1.1) circle [radius=0.15];
\draw [cyan,fill]  (2,0.5) circle [radius=0.15];
\draw [red,fill]  (2,-0.8) circle [radius=0.15];
\draw [red,fill]  (2,-0.5) circle [radius=0.15];
\draw [red,fill]  (2,-1.1) circle [radius=0.15];
\draw [yellow,fill]  (1.7,-0.8) circle [radius=0.15];
\draw [dashed,->] (-0.75,0.25)--(1,0.75);
\draw [dashed,->] (-0.75,-0.25)--(1,-0.75);
\draw [green, fill] (-1.7,-1.4) circle [radius=0.15];
\draw [green, fill] (-1.7,-1.7) circle [radius=0.15];
\draw [green, fill] (2,-1.4) circle [radius=0.15];
\draw [green, fill] (2,-1.7) circle [radius=0.15];
\draw [dashed,->] (-0.75,-1.55)--(1,-1.3);
\end{tikzpicture}
\end{center}
\caption{Pure fragmentation and coagulation processes}\label{Fig1}
\end{figure}
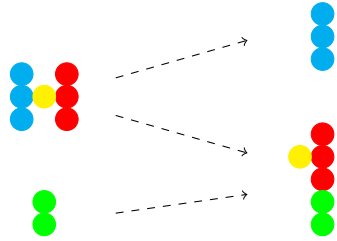

Coagulation and fragmentation equations, modelling the processes of rearrangement of clusters of animate or inanimate matter, are considered to be one of the most fundamental equations of the classical science. These processes consist in splitting large clusters into  smaller blocks and, conversely, forming of larger clusters from smaller ones, see Fig. \ref{Fig1}, and are abundant in natural and technological processes. We can mention here animal groups formation, \cite{DLP17, gueron1996dynamics, okubo1986dynamical}, phytoplankton dynamics, \cite{AcDe03, Ackleh1, BaLa09}, polymer degradation, \cite{montroll1940, simha1941, ZM1, ZiMc85}, planetesimal formation, industrial spray drying, aerosols formation, floculation and many others, see \cite{BLL1} for a more extensive description and references.

The original coagulation--fragmentation equations, formulated by M. Smoluchowski \cite{Smoluch, Smoluch17} for the pure coagulation case and extended to include the fragmentation part in \cite{becker, blatz1945},  assumed the existence of minimum building blocks and thus had the form of an infinite system of ordinary differential equations. However, in many applications such as droplets formation in clouds, fog or aerosols, \cite{Schu40,Scot68}, it makes sense to consider a continuous size variable that can take any positive value. An extension of the original Smoluchowski equations to the continuous case that took the form of an integro--differential equation, was done by M{\"u}ller, see  \cite{muller1928allgemeinen}, and finally Melzak in \cite{Melz57b} incorporated  fragmentation into a continuous particle size model. While the fragmentation term in \cite{becker, blatz1945, Melz57b} is linear, that is, the rate of fragmentation is intrinsically linked to the properties of a given cluster, it is worthwhile to mention the nonlinear fragmentation models, where there breakage is induced by interactions of two clusters, see e.g., \cite{Kost, matv}.

In this paper, we shall deal only with the continuous case and we discard the coagulation processes (but, as explained below, we include transport in the state space). Following original ideas of Smoluchowski, we assume here that a single variable, such as particle mass or size,  is required to differentiate between the reacting particles (a discussion of the cases where other features of the reacting particles, like their shape, play a role in the process can be found in \cite{Watt06}).  Then, for the pure fragmentation process we obtain the following formulation, which can be traced back to the work of  McGrady and Ziff \cite{McZi87} (1987) and Vigil and Ziff \cite{vizi89} (1989),
\begin{numparts}
\label{PhL}
\begin{eqnarray}
u_t(x,t) & = & -a(x) u(x,t) + \  \int_x^{\infty}a(y)b(x,y)
u(y,t)\,d y,  \label{wlcontscfeqn} \\
u(x,0) & = & u_0(x)\ , \qquad x\in \mathbb R_+. \label{PhLicccf}
\end{eqnarray}
\end{numparts}
  Here $u(x,t)$ is the density of particles of size $x > 0$ at time $t$,   $a(x)$ represents the overall rate of fragmentation of an $x$-sized particle, while the coefficient $b(x,y)$, often called the fragmentation kernel, is the distribution function of the
sizes of the daughter particles produced by the parent particle of size $y$. It is assumed to be  nonnegative and measurable, with $b(x,y)=0$ for $x >y$ and
\begin{equation}
\int_0^y xb(x,y)\,
d x = y, \ \mbox{ for each } y > 0,  \label{baleq1}
\end{equation}
 but is otherwise arbitrary.

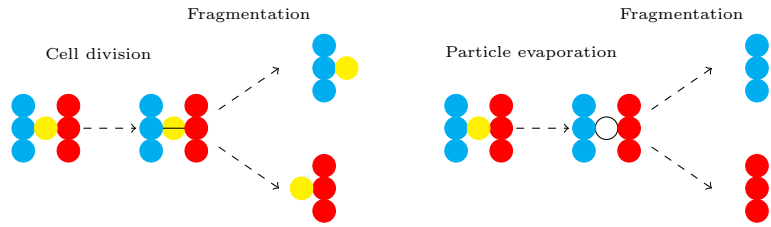
\begin{figure}
\begin{center}
\begin{tikzpicture}[scale=1]
\draw [cyan,fill]  (-2,0) circle [radius=0.15];
\draw [cyan,fill]  (-2,0.3) circle [radius=0.15];
\draw [cyan,fill]  (-2,-0.3) circle [radius=0.15];
\draw [red,fill]  (-1.4,0) circle [radius=0.15];
\draw [red,fill]  (-1.4,0.3) circle [radius=0.15];
\draw [red,fill]  (-1.4,-0.3) circle [radius=0.15];
\draw [yellow,fill]  (-1.7,0) circle [radius=0.15];
\draw [cyan,fill]  (2,0.8) circle [radius=0.15];
\draw [cyan,fill]  (2,1.1) circle [radius=0.15];
\draw [cyan,fill]  (2,0.5) circle [radius=0.15];
\draw [red,fill]  (2,-0.8) circle [radius=0.15];
\draw [red,fill]  (2,-0.5) circle [radius=0.15];
\draw [red,fill]  (2,-1.1) circle [radius=0.15];
\draw [yellow,fill]  (1.7,-0.8) circle [radius=0.15];
\draw[yellow,fill] (2.3,0.8) circle [radius=0.15];
\draw (-1,1) node {\tiny{Cell division}};
\draw (1,1.5) node {\tiny{Fragmentation}};
\draw[yellow,fill] (0,0) circle [radius=0.15];
\draw [-](-0.15,0)--(0.15,0);
\draw [cyan,fill]  (-0.3,0) circle [radius=0.15];
\draw [cyan,fill]  (-0.3,0.3) circle [radius=0.15];
\draw [cyan,fill]  (-0.3,-0.3) circle [radius=0.15];
\draw [red,fill]  (0.3,0) circle [radius=0.15];
\draw [red,fill]  (0.3,-0.3) circle [radius=0.15];
\draw [red,fill]  (0.3,0.3) circle [radius=0.15];
\draw[dashed,->] (-1.2,0)--(-0.5,0);
\draw [dashed,->] (0.6,0.25)--(1.4,0.8);
\draw [dashed,->] (0.6,-0.25)--(1.4,-0.8);
\end{tikzpicture} \qquad \begin{tikzpicture}[scale=1]
\draw [cyan,fill]  (-2,0) circle [radius=0.15];
\draw [cyan,fill]  (-2,0.3) circle [radius=0.15];
\draw [cyan,fill]  (-2,-0.3) circle [radius=0.15];
\draw [red,fill]  (-1.4,0) circle [radius=0.15];
\draw [red,fill]  (-1.4,0.3) circle [radius=0.15];
\draw [red,fill]  (-1.4,-0.3) circle [radius=0.15];
\draw [yellow,fill]  (-1.7,0) circle [radius=0.15];
\draw [cyan,fill]  (2,0.8) circle [radius=0.15];
\draw [cyan,fill]  (2,1.1) circle [radius=0.15];
\draw [cyan,fill]  (2,0.5) circle [radius=0.15];
\draw [red,fill]  (2,-0.8) circle [radius=0.15];
\draw [red,fill]  (2,-0.5) circle [radius=0.15];
\draw [red,fill]  (2,-1.1) circle [radius=0.15];
\draw (-1,1) node {\tiny{Particle evaporation}};
\draw (1,1.5) node {\tiny{Fragmentation}};
\draw[black] (0,0) circle [radius=0.15];
\draw [cyan,fill]  (-0.3,0) circle [radius=0.15];
\draw [cyan,fill]  (-0.3,0.3) circle [radius=0.15];
\draw [cyan,fill]  (-0.3,-0.3) circle [radius=0.15];
\draw [red,fill]  (0.3,0) circle [radius=0.15];
\draw [red,fill]  (0.3,-0.3) circle [radius=0.15];
\draw [red,fill]  (0.3,0.3) circle [radius=0.15];
\draw[dashed,->] (-1.2,0)--(-0.5,0);
\draw [dashed,->] (0.6,0.25)--(1.4,0.8);
\draw [dashed,->] (0.6,-0.25)--(1.4,-0.8);
\end{tikzpicture}
\caption{Fragmentation processes with growth (left) and decay (right)}\label{Fig2}
\end{center}
\end{figure}

Pure fragmentation processes consist only in rearrangements of mass among clusters and thus are (formally) mass conserving. Often, however, the fragmentation is accompanied by events that involve loss or gain of mass, see Fig. \ref{Fig2}. In engineering, physical or chemical practice, mass loss can occur due to oxidation,
melting, sublimation or dissolution of the matter from the exposed surface of the particles. The resulting surface recession
widens the pores of the clusters,  causing their disintegration into  smaller ones, \cite{Ed1, 6, Ed2}.

We may also observe reverse processes, where the clusters, while undergoing fragmentation due to external causes,  may gain mass due to the precipitation of  matter from the environment.  Other examples of this type come from biology, where it has been observed that, depending on circumstances, living organisms form bigger clusters or split into smaller ones, see, e.g., \cite{DLP17, gueron1995dynamics, okubo1986dynamical, Oku, Ackleh1, jackson1990model, RudWiecz2, RudWiecz1}. It is, however, not often fully recognised that the living matter has its own vital dynamics, that is, in addition to forming or breaking clusters, individuals forming them can be born or die inside, which should be adequately represented in the models. In the continuous case, the birth and death processes are incorporated into the model
by adding an appropriate first-order transport term, analogously to the age or size structured McKendrick model, see \cite{Ackleh1, Bana12a, BaLa09,  BaPiRu, PerTr}.  Thus, we shall consider the following modification of (\ref{wlcontscfeqn}) with  either growth ($+$) or decay ($-$),
\begin{numparts}\label{eq3.1}
\begin{eqnarray}
\nonumber
u_t^\pm(x,t) &\pm (r(x)u^\pm(x,t))_x = \\
\label{eq3.1a}
&- a(x)u^\pm(x,t) + \int_x^\infty a(y)b(x,y) u^\pm(y,t)dy,\qquad x,t\in \mathbb{R}_+,\\
\label{eq3.1b}
u^\pm(x,0) &= u_0(x),\qquad x\in \mathbb{R}_+,
\end{eqnarray}
\end{numparts}
where $r$ is a positive and continuous function on $\mathbb R_+$, which gives the rate of growth or decay of clusters in the absence of fragmentation. In the growth case,  if
\begin{equation}
\int_{0}^1 \frac{dx}{r(x)} < +\infty,
\label{rint}
\end{equation}
 then \eref{eq3.1a}--\eref{eq3.1b} must be supplemented by a boundary condition which we assume here to be
\begin{equation}
\lim\limits_{x\to 0^+} r(x)u^+(x) =0.
\label{bc1}
\end{equation}
Problems \eref{wlcontscfeqn}--\eref{PhLicccf} and \eref{eq3.1a}, \eref{eq3.1b} and \eref{bc1} have been the subject of intensive research, both from the practical and theoretical points of view. Nontrivial properties of these problems such as phase transitions (called shattering and runaway fragmentation), see e.g., \cite{Fili61,McZi87, Ed1,6,Ernst}, and the existence of multiple solutions, see e.g., \cite{AizBak,Ba02}, have necessitated the development and application of sophisticated mathematical methods for their analysis. A comprehensive account of the relevant theory can be found in \cite{BLL1}, with further developments given in \cite{BeGa2020,BanLam2020}. In parallel to the theoretical developments, we have seen an extensive search for explicit solutions to \eref{wlcontscfeqn}--\eref{PhLicccf} and \eref{eq3.1a}, \eref{eq3.1b} and \eref{bc1} (as well as to the equations containing the coagulation terms). The interest in explicit solutions stems from several facts -- whenever available, they can be used as benchmark solutions to validate and improve approximate and numerical procedures. They also allow for observing specific phenomena that elude simple theoretical tools and thus stimulate their refinement, and often they give better information about the process that cannot be acquired in a more general context. Thus it is  not surprising that there exists a large body of physical and engineering literature devoted to the derivation and utilization of  explicit solutions. While possibly the first explicit solutions to continuous fragmentation equations are  due to Filippov  \cite{Fili61}, the systematic study of them in the context of polymer degradation is due to Ziff and McGrady, \cite{ZM1, ZiMc85, McZi87}, who derived solutions to the fragmentation equation with the power law coefficients, which are briefly discussed below, in terms of the Kummer confluent hypergeometric function. Their ideas were later extended to some cases of fragmentation equations with decay in a series of papers \cite{6, Ed1, Ed2, huang1996} and used, among others, to determine the shattering regime through the analysis of the scaling properties of the explicit solutions. A related direction of research, based on the knowledge of explicit solutions of fragmentation (and coagulation) equations, deals with finding of the so-called self-similar profiles describing the long term patterns of the evolution of the analysed systems, see \cite{Kapu72a, ZM1, ZiMc85, Pete86, ChRe88, ChRe90, Ziff91, Trea97a, PeRy05, BiTKxx}.

\section{Preliminaries}
Our work, written in the spirit of  the papers \cite{ZM1, ZiMc85, McZi87}, aims at extending their results to a larger class of fragmentation equations with growth or decay. We shall achieve this  by developing a unified approach covering the admissible cases. For the sake of completeness, we begin with the pure fragmentation case.

\subsection{Pure fragmentation models}\label{ss1}

Of special interest in applications are the so-called power-law  coefficients,  defined by
\begin{equation}
a(x) = x^{\alpha}  \mbox{ and } b(x,y) = \frac{\nu + 2}{y}\left(\frac{x}{y}\right)^\nu, \qquad -2 < \nu \leq 0.
\label{pl1}
\end{equation}
For instance, in the context of chain polymers, such coefficients describe the situation in which the scission is
independent of the position along the polymer, but is a function of the length of the chain (if $\alpha\neq 0$), see \cite{ZiMc85}.
Note that the lower limit on $\nu$ is due to the fact that otherwise the integral in  (\ref{baleq1}) is infinite. On the other hand, the upper bound follows from physical features of the model that can be found in \cite[Section 2.2.3.2]{BLL1}. In particular, binary fragmentation  occurs when $\nu = 0$, with the associated initial-value problem  given by
\begin{equation}\label{wlplbinary}
 u_t(x,t) = - x^{\alpha}u(x,t) + 2 \int_x^\infty y^{\alpha-1} u(y,t)\,d y\,,\quad u(x,0) = u_0(x).
\end{equation}
In fact, in \cite{McZi87} it is shown that it is sufficient to investigate this binary equation since
\[
u(x,t) = x^\nu\tilde{u}\left(x^{\frac{\nu + 2}{2}},t\right)
\]
is a solution of the general power-law equation whenever $\tilde{u}$
satisfies the binary equation with $\alpha$ replaced by $2\alpha/(\nu + 2).$   Then, for instance for $\alpha>0$, we can use the integrating factor $e^{-tx^\alpha}$ and the change of variables $\xi = x^\alpha$ so that, letting
\begin{equation}\label{wlu}
u(x,t) = e^{-tx^\alpha}v(\xi,t),
\end{equation}
we can reduce the problem to
\begin{equation}\label{wlredeq}
v_t(\xi,t) = \frac{2}{\alpha}\int_\xi^\infty e^{t(\xi-\eta)}v(\eta,t)\,d \eta, \quad v(\xi,0) = v_0(\xi) = u_0\left(\xi^{\frac{1}{\alpha}}\right).
\end{equation}
 The authors used  a power series technique to obtain a solution, $u_r,$ for the so-called mono-disperse initial state, modelled by  $u_0(x) = \delta_r(x),\ r > 0$, where $\delta$ is the Dirac delta. While such a solution, due to the singularity of the initial condition, is not classical (but can be interpreted in the sense of distributions, see Appendix \ref{secA.2}), it is used as Green's function to obtain the solution to the general initial distribution $u_0$ as
\begin{equation}\label{wlgreenfn}
u(x,t) = \int_{\mathbb{R}_+} u_r(x,t)u_0(r)\,d r.
\end{equation}
 The exact solutions to (\ref{wlplbinary}) obtained in this way were given  in terms of special functions, such as Kummer's confluent hypergeometric function. We shall show that these solutions fit into a more general theory that is presented in this paper.

 \subsection{Fragmentation with growth or decay}\label{ssfgd}

  Here we consider \eref{eq3.1a} and \eref{eq3.1b}, where, in addition to (\ref{pl1}), we  assume
\begin{equation}\label{eq3.1c}
r(x) = k x^\gamma,\qquad a(x) = a x^\alpha,
\end{equation}
with $ k,a>0$, $\alpha\ne0$. In the growth case, if $\gamma>1$, then the characteristics of the problem suffer a finite time blow-up and the total mass of the system becomes infinite in finite time, the pointwise solutions, however,  may still exists, see \cite[Example 5.2.7]{BLL1}. On the other hand, in the decay case the solutions with finite mass exist for any $\gamma\geq0$, \cite[Theorem 9.4]{BaAr}. Further, in the growth case, \eref{rint} is satisfied if $0\leq \gamma <1$ and then we supplement (\ref{eq3.1}) with the boundary condition (\ref{bc1}) (with $r(x) = k x^\gamma$).

First, we simplify the fragmentation part, as in Section \ref{ss1}. Hence, we let
\begin{equation}\label{eq3.2p}
z = a x^\alpha,\qquad v^\pm(z,t) = z^{-\frac{\nu}{\alpha}}u^\pm\left(\left(\frac{z}{a}\right)^{\frac{1}{\alpha}},t\right),\qquad x,t\in\mathbb{R}_+.
\end{equation}
Formal substitution of \eref{eq3.2p} into \eref{eq3.1a}, \eref{eq3.1b} yields
\begin{numparts}\label{eq3.3}
\begin{eqnarray}
\nonumber
v^\pm_t(z,t) \pm \beta z^\mu v^\pm_z(z,t)&=&- \bigl[\pm\theta z^{\mu-1}+  z \bigr]v^\pm(z,t)\\
\label{eq3.3a}
&&
\phantom{xxx}+ m \left\{ \begin{array}{lcl} \int_z^\infty v^\pm(s,t) ds & \textrm{if}& \alpha> 0,\\
\int_0^z v^\pm(s,t) ds & \textrm{if}& \alpha< 0,\end{array}\right.\\
\label{eq3.3b}
&&v^\pm(z,0) = v_0(z) := z^{-\frac{\nu}{\alpha}}u_0\left(\left(\frac{z}{a}\right)^{\frac{1}{\alpha}},0\right),
\end{eqnarray}
for $t, z\in \mathbb{R}_+$, with
\begin{eqnarray}\label{eq3.3c}
\beta &= a^{\frac{1-\gamma}{\alpha}} k\alpha,\qquad \theta = a^{\frac{1-\gamma}{\alpha}}k(\gamma +\nu),\nonumber\\
 m &= \case{(\nu+2)}{|\alpha|}, \qquad \mu = \frac{\gamma +\alpha -1}{\alpha}.
\end{eqnarray}
\end{numparts}
It is important to observe that if $\alpha<0$, then the roles of $+$ and $-$ in (\ref{eq3.3a}) are reversed due to the change of sign of $\beta$ so that the growth problem becomes the decay problem and conversely. This is, however, compensated by the change of the integral operator.

We note that the decay case with $\mu =\theta =0$ was solved in \cite{Ed2, huang1996}. The approach of the authors was to use the characteristic coordinates for the transport part of the equation to remove the derivative with respect to the state variable and work with the simplified equation that resembles the pure fragmentation equation but with time-dependent coefficients, which, in the case with $\mu =\theta =0$, can be incorporated into the unknown function.   Here, we provide a systematic extension of their approach and show that it works for the case of decay and growth with $\mu =1$ and growth with $\mu =0$. We observe that an additional difficulty in the growth case is to make sure that the characteristics fill the whole quadrant $\mathbb R^2_+$ so that the characteristic change of variables is a transformation onto it. Fortunately, in the case of linear growth, this is true,  but in the case of a constant growth rate, we need to impose boundary conditions, which complicates the procedure.

To shorten the notation, we shall define
\[
\mathcal{J}^+[u](x) = \int_x^\infty u(y)dy,\quad
\mathcal{J}^-[u](x) = \int_0^x u(y)dy,\quad x\in\mathbb{R}_+,
\]
and combine the notation into $\mathcal{J}^{\sg}$ where $\textrm{sg}\, \alpha = +$ if $\alpha > 0$ and $\textrm{sg}\, \alpha = -$ if $\alpha<0$.

\begin{figure}
 \begin{center}
 \includegraphics[scale=0.8]{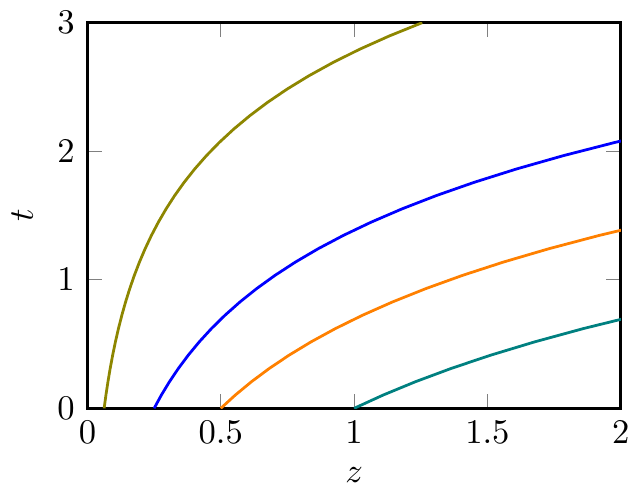} \quad \includegraphics[scale=0.8]{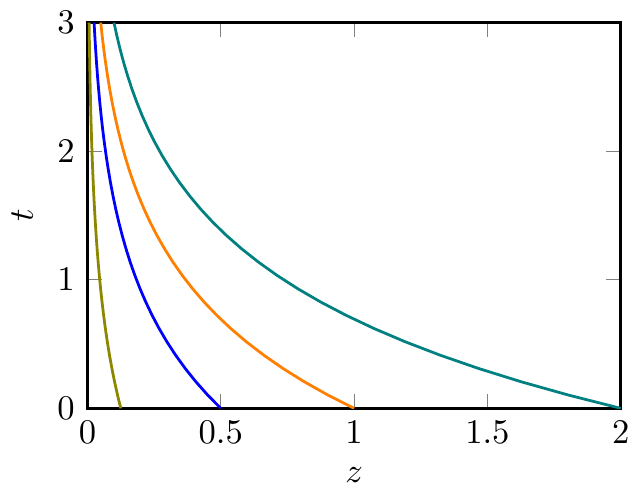}
 \caption{Characteristics determined by the initial condition at $t=0$ fill the first quadrant in the linear growth (left) and decay (right) case.}\label{Fig3}
 \end{center}
 \end{figure}

\paragraph{Linear growth or decay.} Here we assume $\mu =1.$ In terms of the original parameters, this means $\gamma =1,$ so the problems with linear growth/decay become problems with linear growth/decay again (though the growth and decay can switch). Thus,  (\ref{eq3.3a}) becomes
\begin{equation}
v^\pm_t(z,t) \pm \beta z v^\pm_z(z,t)=- \bigl[\pm\theta +  z \bigr]v^\pm(z,t)
+ m \mc J^{\sg} [v^\pm (\cdot,t)](z),
\label{lg1}
\end{equation}
with the same initial condition (\ref{eq3.3b}).
As noted above,  \eref{eq3.1c} implies that in both the growth and the decay scenarios, the characteristics associated
to the transport part of the model fill-up the entire first quadrant $\mathbb{R}_+^2$ so that no boundary
condition at $z=0$ is required, see Fig. \ref{Fig3}.

The characteristics associated to the transport
part of \eref{lg1} are given explicitly by
\[
z^{\pm}(\xi,t) = \xi e^{\pm \beta t},\qquad \xi,t\in\mathbb{R}_+,
\]
where as before the signs $\pm$ correspond to the growth and the decay scenarios in the original variables, respectively.
The substitution
\begin{eqnarray}
&w^\pm(\xi, \tau^\pm) = e^{\tau^\pm\xi} \bigr(1\pm \beta \tau^\pm\bigl)^{\frac{\theta}{\beta}}
v^\pm(z^\pm(\xi,t(\tau^\pm)), t(\tau^\pm)),\nonumber\\
&\tau^\pm(t) = \pm \case{1}{\beta}\bigl(e^{\pm\beta t}-1\bigr),
\quad \tau^\pm\in I^\pm,
\quad I^{\sg} :=\mathbb{R^+},\quad I^{-\sg} := [0,\case{1}{\beta}),\label{trans1}
\end{eqnarray}
(note that an analogous time rescaling was used in \cite{fern}) transforms (\ref{lg1}) to
\begin{numparts}\label{eq3.4}
\begin{eqnarray}
\label{eq3.4a}
& w^\pm_{\tau^\pm}(\xi,\tau^\pm) = m\mathcal{J}^{\sg}
\bigl[e^{-\tau^\pm(\cdot-\xi)} w^\pm(\cdot,\tau^\pm) \bigr](\xi),
\quad (\xi,\tau^\pm)\in\mathbb{R}_+\times I^\pm,\\
\label{eq3.4b}
& w^\pm(\xi,0) = w_0(\xi) :=
\xi^{-\frac{\nu}{\alpha}} u_0\left(\left(\frac{\xi}{a}\right)^{\frac{1}{\alpha}}\right),\qquad \xi\in\mathbb{R}_+.
\end{eqnarray}
\end{numparts}
\paragraph{Constant growth or decay.} Here we assume $\mu =0$ and hence, to avoid singularity at $z=0$, we also require $\theta =0$. In terms of the original parameters, we have $\alpha = 1-\gamma$ and $\gamma = -\nu$.  Thus, the earlier constraints yield $0\leq \gamma <2$ and $-1<\alpha\leq 1$, and hence we can distinguish the following cases:
\begin{enumerate}
\item { Decay with $\alpha<0$.} It requires $\gamma >1$ and no boundary condition in \eref{eq3.1a}--\eref{eq3.1b}; it is transformed into the growth case with $\mc J^-$ requiring boundary condition as $z\to 0^+$ that corresponds to $x\to \infty$. The zero boundary  condition for $v^-$ at $z= 0$ means that $x^\gamma u^-(x)$ converges  0 as $x\to \infty$, in accordance with  \cite[Theorem 9.4]{BaAr}.
    \item {  Growth with $\alpha<0$.} It requires $\gamma >1,$ so no boundary condition in \eref{eq3.1a}--\eref{eq3.1b} is needed; it is transformed into the decay case with $\mc J^-,$ which does not require any boundary condition as $z\to 0^+.$
            \item {  Decay with $\alpha>0$.} It requires $\gamma <1$ and no boundary condition in \eref{eq3.1a}--\eref{eq3.1b}; it is transformed into the decay case with $\mc J^+,$ which does not require any boundary condition as $z\to 0^+.$
                \item {  Growth with $\alpha>0$.} It requires $\gamma <1$ and a boundary condition at $x=0$ in \eref{eq3.1a}--\eref{eq3.1b}; it is transformed into the growth case with $\mc J^+,$ requiring boundary condition as $z\to 0^+$ that corresponds to $x\to 0^+$. The zero boundary  condition for $v^+$ at $z= 0$ means that $x^\gamma u^+(x)$ converges  0 as $x\to 0$.
\end{enumerate}
Having this in mind, the transformed equation takes the form
\begin{equation}
v^\pm_t(z,t) \pm \beta v^\pm_z(z,t)=-  z v^\pm(z,t)
+ m \mc J^{\sg} [v^\pm (\cdot,t)](z).
\label{cg1}
\end{equation}
The initial conditions are again given by (\ref{eq3.3b}), and we will have to supplement (\ref{cg1}) by the boundary condition
\begin{equation}
\lim\limits_{z\to 0^+} v^\pm (z,t) = 0,
\label{cg1bc}
\end{equation}
whenever $\pm \beta>0.$
\begin{figure}
 \begin{center}
 \includegraphics[scale=0.8]{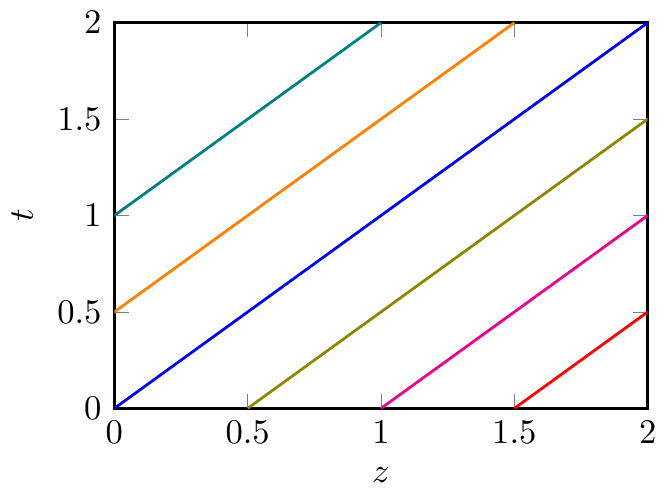} \quad \includegraphics[scale=0.8]{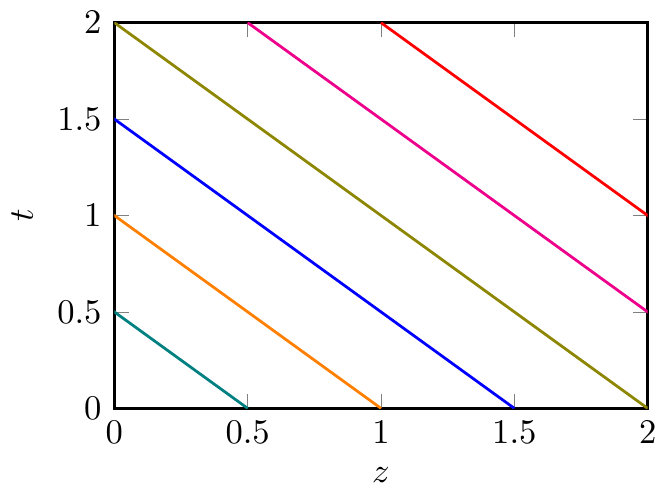}
 \caption{Characteristics determined by the initial condition at $t=0$ fill the first quadrant in the constant decay  case (right). In the growth case (left), the characteristics determined by the initial condition at $t=0$ fill only the region below the main diagonal and hence  a boundary condition at $z=0$ is needed to determine the initial values for the characteristics above it.}\label{Fig4}
 \end{center}
 \end{figure}
Here the characteristics are given by
\begin{equation}\label{charz}
z^\pm(\xi,t) = \pm \beta t + \xi,
\end{equation}
see Fig. \ref{Fig4}.
Then we let $g^\pm(\xi,t)=v(z^\pm (\xi,t),t)$ and hence (\ref{cg1}) becomes
\begin{numparts}\label{G1}
\begin{eqnarray}
&&g^\pm_t(\xi,t)=-(\xi \pm \beta t)g^\pm (\xi, t) + m\mc J^{\sg}[g^\pm(\cdot,t)](\xi),\label{G1a}\\
&&g^\pm(\xi,0)=g_0(\xi) = v_0(\xi),\label{G1b}
\end{eqnarray}
\end{numparts}
for $\xi \in \mathbb R_+$. The way we write the boundary condition depends on the method we choose to solve \eref{G1a}--\eref{G1b}. Here, see Section \ref{sec4}, we extend (\ref{G1}) to $\xi<0$ by assuming that $g_0^\pm(\xi) =0$ for $\xi<0$ and introducing the initial condition $\phi(\xi) = g_0(\xi) + \psi(\xi),$ where $\psi(\xi)=0$ for $\xi >0$ and must be determined for $\xi<0$, see Fig. \ref{Fig5}. In this way, for such an extended solution $g^\pm$, we write (\ref{cg1bc}) as
 \begin{equation}
g^\pm \left(\xi,\mp \beta^{-1}\xi\right ) = 0,\qquad \xi\in \mathbb R_-,
\label{cg1bcg}
\end{equation}
whenever $\pm \beta>0.$ Finally, we set \begin{equation}
w^{\pm}(\xi,t)=e^{\pm\frac{\beta t^2}{2} +\xi t}g^\pm (\xi,t)\label{trans2}\end{equation} and arrive at the equation  analogous to (\ref{wlredeq}) and (\ref{eq3.4a}),
\begin{numparts}\label{G5}
\begin{eqnarray}
\label{G5a}
w^{\pm}_t(\xi, t)&=& m\mc J^{\sg} \bigl[e^{-t(\cdot-\xi)} w^\pm(\cdot,t) \bigr](\xi),\\
\label{G5bp}
w^\pm(\xi,0)&=&\phi(\xi),
\end{eqnarray}
with $\xi >0$ if $\pm\beta <0$ and $\xi\in \mathbb R$ if $\pm\beta >0$ and, in the latter case,
\begin{equation}
w^\pm \left(\xi,\mp \beta^{-1}\xi\right ) = 0,\qquad \xi\in \mathbb R_-.
\label{cg1bcw}
\end{equation}
\end{numparts}
The formal similarity of (\ref{wlredeq}), (\ref{eq3.4a}) and (\ref{G5a}) justifies studying general equations of this type, which is carried out in the following section.

\begin{figure}
\begin{center}
\begin{tikzpicture}[scale=1]
\draw[thin,-] (-3,0)--(3,0);
\draw[thin,-] (0,-1)--(0,0);
\draw (3.2,0) node{\tiny{$\xi$}};
\draw (3.2,0) node{\tiny{$\xi$}};
\draw (0,3.2) node{\tiny{$t$}};
\draw (2,-0.2) node{\tiny{$1$}};
\draw (1,-0.2) node{\tiny{$0.5$}};
\draw (-0.2,-0.2) node{\tiny{$0$}};
\draw (-1,-0.2) node{\tiny{$0.5$}};
\draw (-2,-0.2) node{\tiny{$-1$}};
\draw [purple,-] (2,0)--(2,2.2);
\draw [purple,-] (2,2.6)--(2,3);
\draw [red,-] (1,0)--(1,2.4);
\draw [red,-] (1,2.6)--(1,3);
\draw [blue,-](0,0)--(0,3);
\draw(-0.2,1) node {\tiny{$0.5$}};
\draw(-0.2,2) node {\tiny{$1$}};
\draw [orange,dashed](-1,0)--(-1,1);
\draw [orange,-](-1,1)--(-1,2.2);
\draw [orange,-](-1,2.6)--(-1,3);
\draw [green,dashed](-2,0)--(-2,2);
\draw [green,-](-2,2)--(-2,2.2);
\draw [green,-](-2,2.6)--(-2,3);
\draw [-] (0,0)--(-3,3);
\draw (-3.2,3.2) node {\tiny{$ t= -\xi$}};
\draw (1.5,0.2) node {\tiny{$ v_0(\xi)$}};
\draw (-1.5,0.2) node {\tiny{$ \psi(\xi)$}};
\draw (1.5,2.5) node {\tiny{Region determined}};
\draw (1.51,2.3) node {\tiny{by $ v_0(\xi)$}};
\draw (-1.23,2.5) node {\tiny{Region determined}};
\draw (-1.23,2.3) node {\tiny{by $ w^\pm(\xi,-\xi)=0$}};
\end{tikzpicture}
\caption{Geometry of the problem \eref{G5a}--\eref{cg1bcw} with $\beta =1$ in the growth case in the characteristic coordinates $(\xi,t)$. Recall that $v_0(\xi)$ is only known for $\xi>0$ and $\psi(\xi), \xi<0,$ is to be determined so that $w^\pm \left(\xi,-\xi\right ) = 0$ is satisfied. }\label{Fig5}
\end{center}
\end{figure}
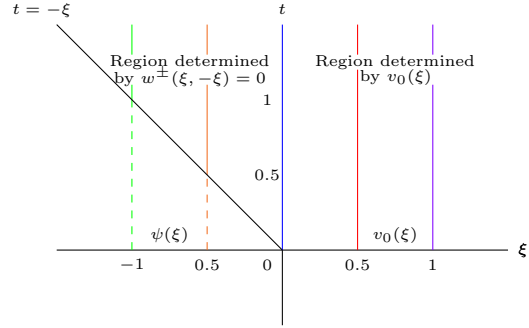

\section{Model equation}\label{sec2}
As we have seen in  Introduction, a large class of fragmentation equations with or without additional
growth or decay can be reduced to the following linear integro-differential equation
\begin{numparts}\label{eq2.1}
\begin{eqnarray}
\label{eq2.1a}
& u_t^\pm(x,t) = m \mathcal{J}^\pm \bigl[ \varphi(\pm t(\cdot-x)) u^\pm(\cdot,t)\bigr](x),
\qquad t\in\mathbb{R}_+,\\
\label{eq2.1b}
& u^\pm(x,0) = u^\pm_0(x),
\end{eqnarray}
with
\[
\mathcal{J}^+[u](x) = \int_x^\infty u(y)dy,\quad
\mathcal{J}^-[u](x) = \int_0^x u(y)dy,
\]
where $x\in\mathbb{R}_+$ or $x \in \mathbb R,$ and the integral kernel is given by an entire function
\begin{equation}\label{eq2.1c}
\varphi(z) = \sum_{n\ge0} \frac{\varphi_n}{n!} z^n,
\end{equation}
\end{numparts}
which we assume to be of finite exponential type $\ell>0,$ see e.g., \cite{Levin1996}. The main aim of this section is
to show that \eref{eq2.1a}--\eref{eq2.1c} is classically solvable for suitably chosen input data $u_0^\pm$.
Our presentation is intentionally abstract. We construct an explicit solution to \eref{eq2.1a}--\eref{eq2.1c}
in a Banach space $X^\pm$, which is chosen so that
\begin{equation}\label{normbound}
\|\mathcal{J}^\pm\|_{X^{\pm}\to X^{\pm}}\le c,
\end{equation}
for some $c>0$.

\paragraph{Remark.} We emphasize that the concrete form of the Banach spaces $X^\pm$ is irrelevant at this point.
All our calculations and formulas hold automatically in every space $X^\pm$, where the basic estimate
\eref{normbound} is satisfied. The concrete choice of $X^\pm$ depends on particular applications, see e.g.
the definition of weighted $L^1$ spaces $X^{\pm\sigma}_{\pm\rho}$ and $X^{+\rho}$ employed in
Sections~\ref{sec3}, \ref{sec4} and Section~\ref{sec4.3}, respectively.

As we shall see shortly, in the functional settings of $X^\pm$ problem \eref{eq2.1a}--\eref{eq2.1c} takes the form
of an abstract linear variable-coefficient ODE.  Indeed, it is easy to verify that
\begin{equation}\label{mcJ}
\left(\mathcal{J}^\pm\right)^n[u](x) = \case{(\pm1)^{n-1}}{(n-1)!}
\mathcal{J}^\pm\bigl[(\cdot- x)^{n-1}u(\cdot)\bigr](x),\qquad n\ge 1.
\end{equation}
Our assumption on the integral kernel $\varphi$ implies that the functions
\begin{equation}
\Phi(z) = \int_{\mathbb{R}_+} \varphi(sz) e^{-s}ds
= \sum_{n\ge 0} \varphi_n z^n,\qquad \widetilde{\Phi}(z) = \sum_{n\ge 0} |\varphi_n| z^n,
\label{Phi}
\end{equation}
are analytic in the disk $\mathbb{D}_{\frac{1}{\ell}} = \bigl\{|z|< \frac{1}{\ell}\bigr\}$.
These observations allow us to define the map
\[
\Phi(z\mathcal{J}^\pm) := \sum_{n\ge 0} \varphi_n \bigr(z\mathcal{J}^\pm\bigl)^n,
\qquad 0\le |z| < \case{1}{c\ell}= :r.
\]
Let $\mathcal{L}(X^\pm)$ denote the space of bounded linear operators on $X^\pm$. By construction,
$\Phi(\pm \cdot \mathcal{J}^\pm): \mathbb{D}_{r}\to \mathcal{L}(X^{\pm })$, is an analytic operator valued
function of $z$ and
\[
\|\Phi(z\mathcal{J}^\pm)\|_{X^{\pm } \to X^{\pm} }\le \widetilde{\Phi}(c|z|).
\]
Hence, locally in time, we may write \eref{eq2.1a}--\eref{eq2.1c} as
\begin{equation}\label{eq2.2}
u_t^\pm = m \Phi(\pm t\mathcal{J}^\pm) \mathcal{J}^\pm u^\pm,\qquad u^\pm(0) = u_0^\pm,\qquad
u^\pm\in C^1((0,r), X^{\pm}).
\end{equation}
The problem \eref{eq2.1a}--\eref{eq2.1c}, written in the abstract form \eref{eq2.2}, is readily integrable.
\begin{lemma}\label{lm2.1}
The problem \eref{eq2.2} is classically solvable. That is, for any $u_0^\pm\in X^{\pm}$,
 there exists a unique classical solution
$u\in C([0,r), X^{\pm})\cap C^1((0,r), X^{\pm})$ satisfying \eref{eq2.2}.
Furthermore, the solutions are given explicitly by
\begin{numparts}\label{eq2.4}
\begin{equation}\label{eq2.3a}
u^\pm(x,t) = \phi_0 u_0^\pm(x) + t \mathcal{J}^\pm[ F(\pm t(\cdot-x)) u_0^\pm(\cdot)](x),
\end{equation}
where the kernel\footnote{Using Hankel's integral representation for the reciprocal of gamma function, one can show that
$F(z) = \int_{\mathbb{R}} \left(\exp\bigl\{m\int_0^{z(1+is)}\Phi(\xi)d\xi \bigr\}
- 1\right) \frac{e^{1+is}ds}{z(1+is)}$.}
$$
F(z) = \sum_{n\ge 0} \case{\phi_{n+1}}{n!(n+1)!} z^n,\qquad
\phi_n = \case{d^n}{dz^n} \exp\Bigl\{m\int_0^z \Phi(\xi)d\xi\Bigr\}\Bigl|_{z=0},
$$
is of finite exponential type $\ell>0$. Equivalently,
\begin{equation}\label{eq2.5}
u^\pm(t) = \phi(t\mathcal{J}^\pm) u_0^\pm,\qquad t\in [0,r),
\end{equation}
where
\begin{equation}
\phi(z) = \exp\Bigl\{m\int_0^z \Phi(\xi)d\xi\Bigr\} = \sum_{n\ge 0} \case{\phi_n}{n!} z^n.
\label{phi}
\end{equation}
\end{numparts}
\end{lemma}
\begin{proof}
(a) The proof is straightforward. First, we rewrite \eref{eq2.2} in the form of Volterra integral equations of the first kind
\begin{equation}\label{eq2.6}
u^\pm(t) = u_0^\pm + m\int_0^t \Phi(\tau\mathcal{J}^\pm) \mathcal{J}^\pm u^\pm(\tau) d\tau,
\qquad 0<t<r.
\end{equation}
Since $\|\Phi(t\mathcal{J}^\pm) \mathcal{J}^\pm\|_{X^{\pm} \to X^{\pm}} \le
c\widetilde{\Phi}(ct)$
and $\widetilde{\Phi}(z)$ is monotone increasing for $t\ge 0$, it follows that the integral equations are uniquely solvable
in $C([0,t_0^\pm], X^{\pm})$, for some $0<t_0^\pm<r$.
Furthermore, since the right-hand side of these equations is in $C^1((0,t^\pm_0), X^{\pm})$,
we conclude that the solutions $u^\pm\in C([0,t_0^\pm], X^{\pm})$
to \eref{eq2.6} are in fact of class $C^1((0,t_0^\pm), X^{\pm})$
and hence satisfy \eref{eq2.2} in the classical sense.

\noindent
(b) To obtain \eref{eq2.3a}, we note that the function $\phi$, defined by (\ref{phi}), is analytic in $\mathbb{D}_{\frac{1}{\ell}}$ as a composition of an entire and an analytic function. Moreover, the reciprocal is given by
$$
\phi^{-1}(z) = \exp\Bigl\{-m\int_0^z \Phi(\xi)d\xi\Bigr\}.
$$
Hence, on account of the commutativity of $t\mc J^\pm$ and $s\mc J^\pm$ for scalar $t,s$,
\[
\phi^{-1}(t\mathcal{J}^\pm)u_t^\pm(t) - \phi^{-1}(t\mathcal{J}^\pm) m\Phi(t\mathcal{J}^\pm)\mathcal{J^\pm} u^\pm(t)
= \frac{d}{dt}\bigl[\phi^{-1}(t\mathcal{J}^\pm)u^\pm(t)\bigr] = 0,
\]
in  $C((0,t_0^\pm), X^{\pm})$. Integrating with respect to time and observing that $\phi^{-1}(t\mathcal{J}^\pm) = \left(\phi(t\mathcal{J}^\pm)\right)^{-1}$,  we infer
$$
u^\pm(t) = \phi(t\mathcal{J}^\pm) u_0^\pm,\qquad t\in [0,t_0^\pm).
$$
In connection with the last formula, we note that the right-hand side of \eref{eq2.5} is analytic when
$t\in \mathbb{D}_{r}$, with values in $X^{\pm}$. Using this fact, it is not difficult
to verify that $u^\pm$, defined in \eref{eq2.5}, satisfies \eref{eq2.2} in the classical sense
for $t\in (0,r)$.
Finally, using the analyticity of $\phi(z)$ and the standard Cauchy estimates, we obtain \eref{eq2.3a}
with $F(z)$ of finite exponential type $\ell>0$.
\end{proof}

Lemma~\ref{lm2.1} is of purely local nature. However, as we shall see in Sections~\ref{sec3} and \ref{sec4},
the integral kernels $F(z)$ appearing in applications to the growth/decay-fragmentation problems are entire
functions of $z$ with moderate (power) growth in $\mathbb{R}_+$. This allows, after some adjustments,
for global extensions of formulae \eref{eq2.3a} and (\ref{eq2.5}) beyond the initial data classes $X^{\pm}$
and the time interval $(0,r)$.

\section{Solution to the fragmentation models with linear growth and decay rates}\label{sec3}

Since it is obvious that $\varphi(z)  := e^{- (\sg) z}$ is of finite exponential type $1$,
the model \eref{eq3.4a}--\eref{eq3.4b} falls in the scope of the theory presented in Section~\ref{sec2}. In particular,
if we define
\begin{equation}\label{spaces1}
X^{\pm\sigma}_{\pm \rho} : = L^1(\mathbb{R}_+,x^{\pm\sigma} e^{\pm\rho x}dx), \quad
\|u \|_{X^{\pm \sigma}_{\pm \rho}}
= \int_{\mathbb{R}_+}\!\!\! |u|(x)x^{\pm \sigma} e^{\pm \rho x} dx,
\end{equation}
for $\sigma\geq 0,\rho> 0,$ then direct calculations show that
\[
\|\mathcal{J}^\pm\|_{X^{\pm \sigma}_{\pm \rho}\to X^{\pm\sigma}_{\pm\rho}}\le \case{1}{\rho},
\qquad \rho>0,\qquad \sigma\ge 0.
\]
We note that the spaces with exponential weights for the solution of \eref{eq3.4a}--\eref{eq3.4b} are natural as the later problem is obtained from the original ones through exponential scalings such as (\ref{wlu}), (\ref{trans1}) or (\ref{trans2}). The solutions to the original problems are considered, however, in spaces with no exponential weights, as seen in Section \ref{sec3.2}.

It immediately follows from
Lemma~\ref{lm2.1} that the local in time classical solutions to \eref{eq3.4a}--\eref{eq3.4b}, with the initial data
$w_0^\pm\in X_{\sg\,\rho}^{\sg\,\sigma}$, $\rho>0$, $\sigma\ge 0$, are given explicitly by
\begin{eqnarray}
\nonumber
w^\pm(\xi,\tau^\pm) &= w_0^\pm(\xi) \\
\label{eq3.2}
&+ m\tau^\pm \mathcal{J}^{\sg} \Bigl[
{}_1F_1\Bigl(1-(\sg) {m};2;- \tau^\pm(\cdot-\xi) \Big) w_0^\pm\Bigr](\xi),
\end{eqnarray}
where
\[
{}_1F_1(a;b;z) = \sum_{n\ge 0} \case{(a)_nz^n}{(b)_nn!},
\]
is the Kummer confluent hypergeometric function of the first kind,
\cite[Formula 13.1.2, p.~504]{abramowitz1964handbook}.
We note also that in view of the Kummer transformation, see \cite[Formula 13.1.27, p.~505]{abramowitz1964handbook},
\begin{equation}\label{eq3.3p}
{}_1F_1(a;b;z) = e^z{}_1F_1(b-a;b;-z),
\end{equation}
$w^\pm$, defined in \eref{eq3.2}, are nonnegative for nonnegative input data $w_0^\pm$.

\paragraph{Remark.}
Alternatively, to use \eref{eq2.5}, we observe that, since $\varphi(z) = e^{-(\sg)z}$,
$$
\Phi(z) = \frac{1}{1+(\sg)z}.
$$
While the series defining $\Phi$ in (\ref{Phi}) converges only for $|z|<1$, $\Phi$ is analytic everywhere except for $z=-(\sg)1$.
Thus, $\phi(z) = (1+(\sg)z)^{(\sg)m}$ where, for non-integer $m>0$ (as well as for the function
$\ln (1+(\sg)z)$ that appears in intermediate calculations), we assume that the respective functions are defined in
$\mathbb C$, cut along the ray $(-\infty, -1]$, when $\alpha>0$, or along $[1,+\infty)$, when $\alpha<0$.
Then, as long as the spectrum of $\tau^\pm \mc J^{\sg}$ does not intersect the respective line, we have
\begin{equation}\label{eq3.4abstract}
w^\pm (\xi, \tau^\pm) = (I +(\sg)\tau^\pm \mc J^{\sg} )^{(\sg)m}[w^\pm_0](\xi).
\end{equation}
In the particular case of $\alpha>0$, this approach yields global in $\tau^\pm\in\mathbb{R}_+$ solutions $w^\pm$.
The rigorous proof of this fact is given in Appendix~\ref{appB}, Proposition~\ref{propB1}.

\subsection{Explicit solutions}\label{sec3.1}
Backward substitution shows that the solutions to \eref{eq3.1a}--\eref{eq3.1b} are given explicitly by
\begin{eqnarray}
\nonumber
u^{\pm}(x,t) &= \exp\Bigl\{\mp kt\mp\case{a x^\alpha}{k\alpha}(1-e^{\mp k\alpha t}) \Bigr\}
\biggl[ u_0^\pm(x e^{\mp kt})\\
\nonumber
&\pm \case{a(\nu+2)}{k\alpha}(1-e^{\mp k\alpha t})
\int_x^\infty
{}_1F_1\Bigl(1-\case{\nu+2}{\alpha};2;\case{\mp a}{k\alpha}(1-e^{\mp k\alpha t})(y^\alpha-x^\alpha)\Bigr)\\
\label{eq3.5a}
&x^\nu y^{\alpha-\nu-1} u^\pm_0(y e^{\mp kt})dy\biggr].
\end{eqnarray}
We remark at this point that \eref{eq3.5a} follows from \eref{eq3.2} in a purely formal manner, hence it
requires proper justification. To avoid overloading text with huge amount of mathematical technicalities,
the rigorous proofs are postponed to Appendix~\ref{App1}. Here we mention that \eref{eq3.5a}
are indeed genuine solutions in the sense of distributions. That is, these solutions satisfy \eref{eq3.1a}--\eref{eq3.1b}
in the space of Schwartz distributions $\mathcal{D}'(\mathbb{R}^2_+)$, see \cite{Bre}, for initial data $u_0^\pm$ in $\mathcal{D}'(\mathbb{R}_+)$.
As an immediate consequence, for the monodisperse initial
data $u^\pm_0 = \delta_{x_0}$, $x_0\in\mathbb{R}_+$,
we have\footnote{In formula \eref{eq3.5b} and everywhere
below, $\chi_{A}(\cdot)$ denotes the characteristic function of set $A\subset\mathbb{R}_+$.}
\begin{eqnarray}
\nonumber
u^\pm(x,t) &= \exp\Bigl\{\mp\case{a x_0^\alpha}{k\alpha}(e^{\pm k\alpha t}-1) \Bigr\}
\biggl[\delta_{x_0e^{\pm kt}}(x)\\
\nonumber
& \pm \chi_{[0,x_0e^{\pm kt}]}(x) \case{a(\nu+2)}{k\alpha}(e^{\pm k\alpha t}-1)x_0^{\alpha-\nu-1}x^\nu \\
\label{eq3.5b}
&\;\;
{}_1F_1\Bigl(1+\case{\nu+2}{\alpha};2;\case{\mp a}{k\alpha}(1-e^{\mp k\alpha t})
(x^\alpha-x_0^\alpha e^{\pm k\alpha t})\Bigr)\biggr].
\end{eqnarray}
Formula \eref{eq3.5b} has straightforward physical interpretation. Expanding the right-hand side
of \eref{eq3.5b}, we see that the term containing $\delta$-distribution describes the evolution and mass loss
of the original particle of size $x_0$, while the term with the Kummer function provides a continuous mass
distribution of daughter particles generated by multiple fragmentation and transport events.

\paragraph{Remark.} As usual, stronger properties of solutions are obtained for regular initial data $u_0^\pm$.
For instance, if $(u^\pm_0)_x\in X_0^{p+1}$, $u^\pm_0\in X_0^{p+\alpha}$ and
\begin{numparts}\label{eq3.18}
\begin{eqnarray}
\label{eq3.18a}
&\alpha>0,\qquad p\ge \alpha-\nu-1,\\
\label{eq3.18b}
&\alpha<0,\qquad p\ge 0,\qquad p>1+\alpha,
\end{eqnarray}
\end{numparts}
then
\begin{eqnarray}
\label{eq3.16bp}
&u^\pm_t, (ru^\pm)_x, au \in C\bigl((0,T),X^{p}_0\bigr),
\end{eqnarray}
and \eref{eq3.5a} satisfies \eref{eq3.1a}--\eref{eq3.1b} in the classical sense of $X_0^p$.
For the proof of this and related facts the reader is referred to Appendix \ref{App1}, Corollary~\ref{lm3.2}.

\subsection {Moments}\label{sec3.2}

Moments $M_p^\pm(t) := \|u^\pm(\cdot,t)\|_{X^p_0}$
of the solutions are of physical importance as they provide information about the global state of the evolving system.
For instance, the zeroth moment $M_0^\pm(t)$ gives the number of particles in the system at time $t$,
the first moment $M_1^\pm(t)$ describes the evolution of the total mass of the system, while the higher order
moments are related to the distribution of mass between small and large clusters. Furthermore, the behaviour
of the first moment is related to the occurrence of a phase transition phenomenon,
known as "shattering", \cite{Fili61,McZi87}, that describes an unaccounted for loss of mass from the system; we mentioned it in Introduction.
In the context of fragmentation with growth or decay, shattering refers to the fact that  the evolution of the total mass
of the system is not determined solely by the mass growth/decay terms $\pm (r(x)u^\pm(x,t))_x$ built into the model,
\cite{Ed1,6} and \cite[Sections 5.2.7\&8]{BLL1}.

As mentioned in Section~\ref{sec3.1}, for integrable input data $u_0^\pm\in X_0^p$, the nonnegative $p$-th
order moments $M_p^\pm(t)$ are well defined and remain finite at each instance of time $t\ge 0$, only
if the inequalities \eref{eq3.18a}-\eref{eq3.18b} are satisfied. There is, however, an interesting difference in the
behaviour of higher and lower order moments of $u^\pm$, $t>0$, when $\alpha>0$ and when $\alpha<0$.
It follows from \eref{eq3.16a} of Corollary~\ref{lm3.2} that in the first case all higher order moments
$p\le q\le p+\nu+2$ become finite instantly at $t>0$. In contrast, for $\alpha<0$ only the lower order
moments $1+\alpha<q\le p$ and $0\le q$, remain well defined. This phenomenon is related to
two different types of moment regularization effect (discovered recently in \cite{BeGa2020, BanLam2020}) induced
by the fragmentation rate $a$, with $\alpha>0$ and $\alpha<0$.

\begin{figure}[t!]
\begin{center}
\includegraphics{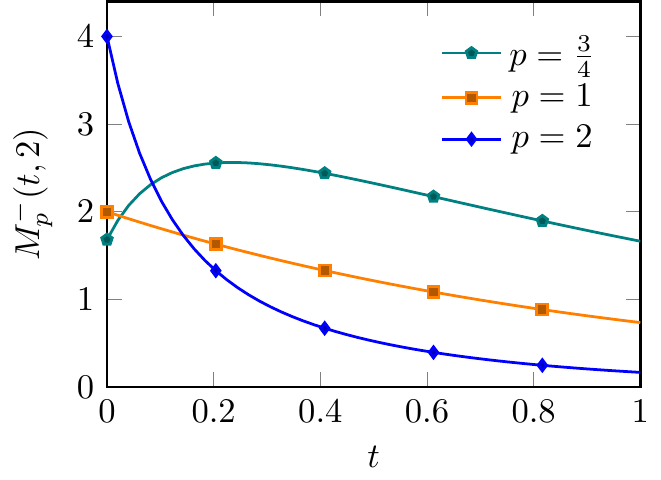} \hspace{2mm} \includegraphics{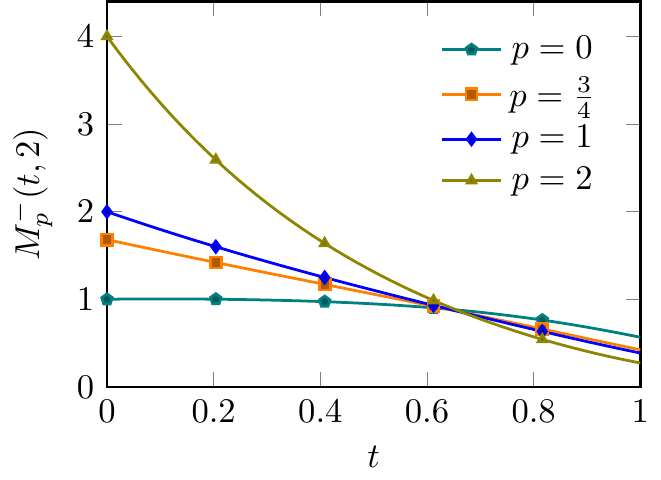}\\
\includegraphics{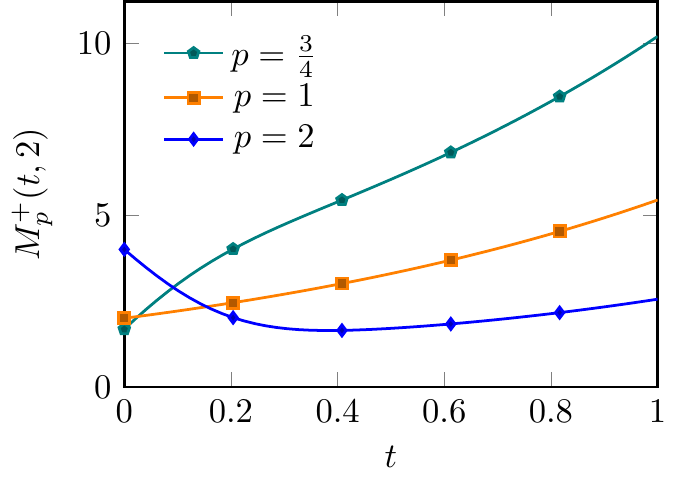} \hspace{2mm} \includegraphics{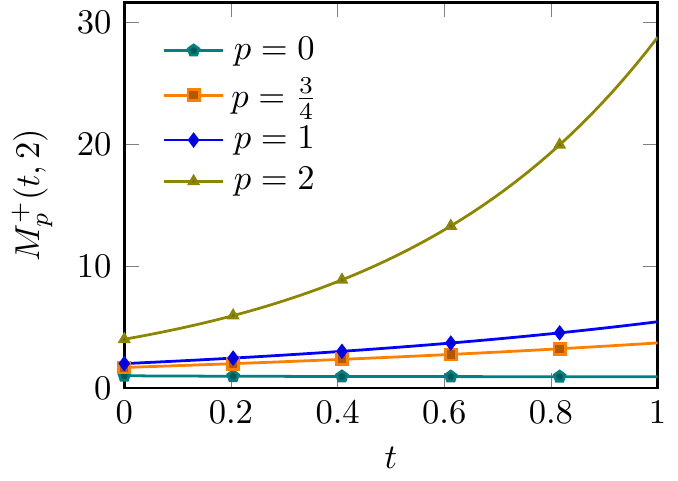}
\end{center}
\caption{Evolution of moments with $u^\pm_0(x) = \delta_2(x)$ and $k=1$ $a=1$, $\nu=-\frac{3}{2}$:
left column $\alpha=3$; right column $\alpha=-3$.}\label{fig1}
\end{figure}

\paragraph{Remark.}
The explicit formulae \eref{eq3.5a} and \eref{eq3.5b} allow for direct computations of nonnegative-order
moments for integrable and monodisperse input data $u_0^\pm$, respectively. However, in view of the linearity
of the model \eref{eq3.1a}--\eref{eq3.1b} and of the moment functionals $M_p^\pm(t)$, $p\ge 0$,
we present calculations only for the latter case. For general integrable data $u^\pm_0\in X_0^p$, the
dynamics of $M_p(t)$ can be read off the monodisperse case via the standard superposition principle
(see e.g., \eref{wlgreenfn}).

To emphasize the dependence on the initial data, in the monodisperse case we denote
$M_p^\pm(t,x_0) := \|u^\pm(\cdot,t)\|_{X_0^p}$. Then, integrating \eref{eq3.5b} with the aid of identities
\eref{eq3.10b} and \eref{eq3.11}, for the general value of $p\ge 0$, satisfying \eref{eq3.18a}--\eref{eq3.18b},
we get
\begin{numparts}\label{eq3.19}
\begin{eqnarray}
\nonumber
M_p^\pm(t,x_0) &
=& \exp\bigl\{\pm pkt\mp\case{a x_0^\alpha}{k\alpha}(e^{\pm k\alpha t}-1)\bigr\}
x_0^p\\
\label{eq3.19a}
&&\cdot{}_1F_1\Bigl(\case{\nu+2}{\alpha};\case{p+\nu+1}{\alpha};
\case{\pm a x_0^\alpha}{k\alpha}(e^{\pm k\alpha t}-1) \Bigr),\qquad\alpha>0,
\end{eqnarray}
and
\begin{eqnarray}
\nonumber
M_p^\pm(t,x_0) &
=& \frac{\Gamma\bigl(\frac{\alpha-p+1}{\alpha}\bigr)}{\Gamma\bigl(\frac{\alpha-p-\nu-1}{\alpha}\bigr)}
\exp\bigl\{\pm pkt \mp\case{a x_0^\alpha}{k\alpha}(e^{\pm k\alpha t}-1)\bigr\}x_0^p\\
\label{eq3.19b}
&&\cdot\Psi\Bigl(\case{\nu+2}{\alpha};\case{p+\nu+1}{\alpha};
\case{\pm a x_0^\alpha}{k\alpha}(e^{\pm k\alpha t}-1) \Bigr),\qquad\alpha<0,
\end{eqnarray}
\end{numparts}
where $\Psi(a;b;z)$ is the Kummer\footnote{In fact, $\Psi(a;b;z)$ was introduced by F. Tricomi and
accordingly, in some texts, is called the Tricomi hypergeometric function.}
hypergeometric function of the second kind, see \cite[Formulas 13.2.5 and 13.1.29, p.~505]{abramowitz1964handbook}. In particular, letting $p=1$ in \eref{eq3.19a}-\eref{eq3.19b}, we see that the total mass of the particle system associated to a monodisperse initial distribution evolves according to the formulae
\begin{numparts}\label{eq3.20}
\begin{eqnarray}
\label{eq3.20a}
&M_1^\pm(t,x_0) = e^{\pm kt}x_0,\qquad\alpha>0,\\
\label{eq3.20b}
&M_1^\pm(t,x_0)
= \frac{e^{\pm kt}x_0}{\Gamma\bigl(1-\frac{\nu+2}{\alpha}\bigr)}
\Gamma\Bigr(1-\case{\nu+2}{\alpha};
\case{\pm a x_0^\alpha}{k\alpha}(e^{\pm k\alpha t}-1) \Bigr),\qquad\alpha<0,
\end{eqnarray}
\end{numparts}
where $\Gamma(a;z) = \int_z^\infty e^{-s} s^{a-1} ds,\,a>0,$ is the incomplete gamma function (see \cite[Formula 6.5.3, p.~260]{abramowitz1964handbook}).
Since  \eref{eq3.20a} describes also the evolution of the total mass due to the transport terms $\pm k(xu(x,t))_x$,
we immediately see that there is no shattering if $\alpha>0$, while \eref{eq3.20b} shows that shattering occurs
in both decay and the growth scenarios for $\alpha<0$.
However, in the growth case, on account of the asymptotic identity
\[
\Gamma(a;x) = \mathcal{O}(x^{a-1}e^{-x}),\quad x\to+\infty,
\]
\cite[Formula 6.5.32, p.~263]{abramowitz1964handbook},  for large values of $t$  shattering
is dominated by the linear growth and in this case we have
\[
\lim_{t\to\infty} M^+_1(t,x_0)=\infty,\qquad x_0\in\mathbb{R}_+,\qquad \alpha<0.
\]
The typical behavior of moments in all four cases covered by \eref{eq3.1a}--\eref{eq3.1b}, with the monodisperse
initial data $u^\pm_0(x)=\delta_{2}(x)$, is shown in Fig~\ref{fig1}.

\subsection{Non-uniqueness}\label{sec3.3}
It is well known that pure fragmentation ($k=0$) equation \eref{eq3.1a}--\eref{eq3.1b}, with $\alpha>0$,
has multiple solutions satisfying the same initial data. This phenomenon was first observed in \cite{AizBak}, and it was explained in \cite{Ba02, Bana04}, where the author also shown that it is related to the non-maximality of the generator of the semigroup associated with the fragmentation equation.

It turns out that similar situation occurs in the context
of growth/decay-fragmentation model \eref{eq3.1a}-\eref{eq3.1b}, for general $k>0$. Indeed, separating
variables in \eref{eq3.4a} and using transformation \eref{eq3.17a}-\eref{eq3.17b}, we infer that the family
of functions
\begin{equation}\label{eq3.21}
\hat{u}^\pm(x,t) = x^\nu\int_{0}^\infty\frac{
\exp\bigl\{\mp k(\nu+1)t \pm \case{\mu}{k\alpha}(e^{\pm k\alpha t}-1)\bigr\}}{{ \bigl(\frac{\mu}{a}+x^\alpha
e^{\mp k\alpha t}\bigr)^{\frac{\alpha+\nu+2}{\alpha}}}}\hat{u}_0^\pm(\mu) d\mu,
\end{equation}
with $ \alpha>0$, satisfy \eref{eq3.1a} pointwise, for any zeroth order distribution $\hat{u}^\pm_0\in \mathcal{D}'(\mathbb{R}_+)$.
It is a trivial exercise
to verify that these solutions are $p$-integrable, provided $-(1+\nu)< p<1+\alpha$; and are
classical in $X_0^p$, provided $-(1+\nu)<p<1$. In the former case, the $p$-th order moments
$\hat{M}_p^\pm(t) :=\|\hat{u}^\pm(\cdot,t)\|_{X_0^p}$ are well defined
and are given by the formula
$$
\hat{M}_p^\pm(t) =
\case{a^{\frac{2-p}{\alpha}}}{\alpha}B\bigl(\case{p+\nu+1}{\alpha},\case{\alpha-p+1}{\alpha}\bigr)\int_{0}^\infty
\exp\bigl\{\pm pkt \pm \case{\mu}{k\alpha}(e^{\pm k\alpha t}-1)\bigr\}\mu^{\frac{p-2}{\alpha}}
\hat{u}_0^\pm(\mu)d\mu.
$$
Hence, we see that in both the growth and the decay scenarios, the total mass of the
system described by \eref{eq3.21} is amplified by the spurious factor
$\exp\bigl\{\pm \case{\mu}{k\alpha}(e^{\pm k\alpha t}-1)\bigr\}$, rendering these solutions
physically infeasible.

\section{Solution to fragmentation models with constant growth and decay rates}\label{sec4}

Since the solutions here are also based on Lemma \ref{lm2.1} and the coordinate transformations defined by
the characteristics are diffeomorphisms, as in \eref{eq3.17a}--\eref{eq3.17b}, in this section we shall focus only
on deriving the formulae for solutions to \eref{G5a}--\eref{cg1bcw}.

We have four different cases described in paragraph \textit{Constant growth or decay} of Section \ref{ssfgd}.
However, the cases (ii) and (iii) do not require boundary conditions and both are confined to the first quadrant,
as in Section \ref{sec3}.

\subsection{Cases $\pm\beta <0$.}\label{sec4.1}
This scenario covers items (ii) and (iii). That is, the decay case with $\alpha>0$ or the growth case with $\alpha<0$
in \eref{eq3.1a}--\eref{eq3.1c}.  Here, problem \eref{G5a}--\eref{G5bp} takes the form
\begin{numparts}\label{G55}
\begin{eqnarray}
w^{\pm}_t(\xi, t)&=& m\mc J^{\sg} \bigl[e^{-t(\cdot-\xi)} w^\pm(\cdot,t) \bigr](\xi),\label{G5aa}\\
w^\pm(\xi,0)&=&w_0(\xi),\qquad \xi\in \mathbb R_+.\label{G5ab}
\end{eqnarray}
\end{numparts}
In this case, analogously to \eref{eq3.4abstract}, the local solutions are given by
\begin{equation}
w^{\pm}(\xi, t) = (I+(\sg)t\mc J^{\sg})^{(\sg)m} [w_0](\xi).
\label{solw11}
\end{equation}
Situation here is completely identical to that considered in detail in Section~\ref{sec3}. For that reason,
our presentation here is very laconic. We mention that in view of Lemma \ref{lm2.1} the solutions \eref{solw11}
are classical for small values of $t>0$, in the sense of spaces $X^{\sg \sigma}_{\sg \rho}$ from \eref{spaces1}.
Further, on the account of Proposition~\ref{propB1}, \eref{solw11}, with $\alpha>0$, holds for any finite
value of $t>0$. The explicit solutions to \eref{eq3.1a}--\eref{eq3.1c} follow from \eref{solw11}
by passing back from the characteristic $(\xi,t)$ back to the physical $(x,t)$ variables
(for the concrete formulae see Section~\ref{sec4.4} below). Repeating almost verbatim the extension arguments of
Lemma~\ref{lm3.1}, it is not difficult to verify that for the resulting explicit
formulae all conclusions of Corollary~\ref{lm3.2} hold.

\subsection{The case $\pm\beta >0$ and $\alpha<0$.}\label{sec4.2}
This is only possible if $\beta<0$, that is, we deal with (i) -- the decay with $\alpha<0$ in
\eref{eq3.1a}--\eref{eq3.1c}. Then problem \eref{G5a}--\eref{G5bp} takes the form
\begin{numparts}\label{G5b}
\begin{eqnarray}
w^{-}_t(\xi, t)&=& m\mc J^{-} \bigl[e^{-t(\cdot-\xi)} w^-(\cdot,t) \bigr](\xi)
= m\int_0^\xi e^{-t(\eta-\xi)}w(\eta,t)d\eta,\label{G5ba}\\
w^-(\xi,0)&=&\phi(\xi), \label{G5bb}
\end{eqnarray}
with $\xi \in \mathbb R,$ $\phi(\xi) = v_0(\xi) + \psi(\xi),$ where $v_0$ is assumed to be extended by 0 to
$\mathbb R_-,$   $\psi(\xi)=0$ for $\mathbb R_+$ and must be determined  on $\mathbb R_-$, so that
\begin{equation}
w^- \left(\xi,\beta^{-1}\xi\right ) = 0,\qquad \xi\in \mathbb R_-.
\label{cg1bcwa}
\end{equation}
\end{numparts}
Having in mind the equivalent representation of the solution, given by \eref{eq2.3a}, we see that
$(I-t\mc J^-)^{-m}[v_0](\xi)$ converges to zero as $\xi\to 0^+$ if so does $v_0$
(which, due to $\alpha<0$, is equivalent to $u_0$ vanishing as $x\to \infty$). This can be ascertained,
at least for small $t>0$, by taking the series expansion of $(I-t\mc J^-)^{-m}[v_0](\xi)$ and noting that
its terms are $v_0$ and integrals from $0$ to $\xi$; see also \eref{BRP} for integer values of $m$.
Hence, $(I-t\mc J^-)^{-m}[v_0](\xi)$ can be continuously extended by $0$ to $\mathbb R_-$ for any $t\geq 0$
and a quick reflection leads us to the conclusion that the function
\begin{equation}\label{solwm}
w^-(\xi,t) = \left\{\begin{array}{lcl} 0 &\mathrm{for}& \xi\in \mathbb R_-, \\
(I-t\mc J^-)^{-m}[v_0](\xi) &\mathrm{for}& \xi\in \mathbb R_+,
\end{array}\right.
\end{equation}
 is a solution to \eref{G5ba}.

\paragraph{Remark.} This result is natural if we take into account \cite[Theorem 9.4]{BaAr} and comment (i) in
Section \ref{ssfgd}. Indeed, according to the former, if $\gamma >1$, then the characteristics of the transport
part fill only the region in $\mathbb R_+\times \mathbb R_+$ bounded  by the limit characteristic
$t= \frac{x^{1-\gamma}}{k(\gamma-1)}$ and the solutions vanish identically outside of it.
Since for $\mu =0$, we have $\alpha = 1-\gamma$, on the account of \eref{eq3.2p} and \eref{eq3.3c},
this characteristic is transformed into
$$
t = -\frac{1}{k\alpha a}ax^\alpha = -\frac{1}{\beta} z^-,
$$
which is precisely the limiting characteristic \eref{charz}, separating the region of influence of the initial condition
from the region of influence of the boundary condition in \eref{cg1bc}.

\paragraph{Example.} In this example, we shall present the application of \eref{solwm} to a particularly simple
case of \eref{eq3.1a}--\eref{eq3.1c}, illustrating the pertinent techniques.
Consider
\begin{numparts}\label{eq5.1}
\begin{eqnarray}
u_t^-(x,t) &- (x^{\frac{4}{3}}u^-(x,t))_x = \label{eq5.1a}
- x^{-\frac{1}{3}}u^- + \frac{2}{3}x^{-\frac{4}{3}}\int_x^\infty u^-(y,t)dy,\\
\label{eq5.1b}
u^-(x,0) &= \delta_{x_0}(x),
\end{eqnarray}
\end{numparts}
where $x_0,x,t\in\mathbb{R}_+$. Here, $\gamma = \frac{4}{3} = -\nu$ and $\alpha=-\frac{1}{3}$,
leading to $m=2$, $\beta = -\frac{1}{3}$ so that the original characteristics
$\frac{1}{\sqrt[3]{x}}-\frac{1}{\sqrt[3]{\eta}} = 3t$ are transformed into $z - \frac{1}{3}t = \xi$ and
the limiting characteristic $\frac{1}{\sqrt[3]{x}}= 3t$ (as $\eta \to \infty$ transforms into $z= \frac{1}{3}t$
(with $\xi =0$). Moreover, using \eref{eq3.2p}, the initial condition for $u^-$ is transformed into
$v_0(\xi)= \frac{1}{3}\delta_{x_0^{-1/3}}(\xi) = \frac{1}{3}\delta(\xi-x_0^{-1/3})$. Precisely, by \eref{eq3.3c},
for any test function $\phi$, we have
\begin{eqnarray*}
\int_0^\infty v_0(z)\phi(z) dz &=& \int_0^\infty z^{-4} \delta(z^{-3}-x_0) \phi(z) dz
= \frac{1}{3}\int_{0}^\infty \delta(r-x_0) \phi(r^{-\frac{1}{3}}) dr\nonumber\\
& =& \frac{1}{3}\phi(x_0^{-\frac{1}{3}})
= \frac{1}{3}\int_0^\infty \delta(r-x_0^{-\frac{1}{3}}) \phi(r) dr.\nonumber
\end{eqnarray*}
Now, using \eref{res1} with $\lambda =1$, we find that
\begin{eqnarray}
&(I-t\mc J^-)^{-m}[g](\xi) = g(\xi)
+ \sum_{n=1}^m\binom{m}{n} t^n \left(\mc J^-\right)^n[e^{-t(\cdot-\xi)}g(\cdot)](\xi)\label{BRP}\\
&\qquad\qquad
 = g(\xi)
+ \sum_{n=1}^m \binom{m}{n} \frac{(-1)^{n-1} t^n}{(n-1)!}
\int_0^\xi e^{-t(\sigma-\xi)}(\sigma-\xi)^{n-1}g(\sigma)d\sigma,
\nonumber
\end{eqnarray}
where \eref{mcJ} was employed. Hence, using \eref{solwm} with \eref{eq5.1b} and \eref{BRP}
with $m=2$, for $\xi>0$,
$$
w^-(\xi,t) = \frac{1}{3}\delta\left(\xi -\frac{1}{\sqrt[3]{x_0}}\right) +
 \frac{te^{t\left(\xi-\frac{1}{\sqrt[3]{x_0}}\right)}}{3}
\left(2+t\left(\xi-\frac{1}{\sqrt[3]{x_0}}\right)\right)\chi_{[x_0^{-1/3},\infty)}(\xi),
$$
where $\chi$ is, as before, the characteristic function of the indicated set. Using again
the rules of transformation of the delta function, we arrive at the solution
\begin{eqnarray*}
u^-(x,t) &= &e^{-t\left(\frac{1}{\sqrt[3]{x_0}}+\frac{t}{6}\right)}\left(\delta\left(x - \frac{27 x_0}{(tx_0^{1/3} +3)^3}\right)\right. \\
&+&\left.
\frac{t}{3x^{\frac{4}{3}}}
\left(2+t\left(\frac{1}{\sqrt[3]{x}}-\frac{1}{\sqrt[3]{x_0}}-\frac{t}{3}\right)\right)
\chi_{\left[0, \frac{27 x_0}{(tx_0^{1/3} +3)^3}\right]}(x)\right).
\end{eqnarray*}

\subsection{The case $\pm\beta >0$ and $\alpha>0$.}\label{sec4.3}
 In this scenario necessarily $\beta>0$ and we have the growth case with $\alpha>0$ in \eref{eq3.1a}, that is,
case (iv). Then the problem \eref{G5a}--\eref{G5bp} takes the form
 \begin{numparts}\label{G5c}
\begin{eqnarray}
w^{+}_t(\xi, t)&=& m\mc J^{+} \bigl[e^{-t(\cdot-\xi)} w^+(\cdot,t) \bigr](\xi)
= m\int_\xi^\infty e^{-t(\eta-\xi)}w(\eta,t)d\eta,\label{G5ca}\\
w^+(\xi,0)&=&\phi(\xi),
\label{G5cb}
\end{eqnarray}
with $\xi \in \mathbb R,$ $\phi(\xi) = v_0(\xi) + \psi(\xi),$ where, as before, $v_0$ is assumed to be extended by $0$
to $\xi<0,$   $\psi(\xi)=0$ for $\xi >0$, but must be determined in a nontrivial way for $\xi<0$, so that
\begin{equation}
w^+ \left(\xi,-\beta^{-1}\xi\right ) = 0,\qquad \xi<0,
\label{cg1bcwb}
\end{equation}
\end{numparts}
is satisfied. It is not difficult to verify that the operator $\mc J^+$ is bounded on the space
$X^{+\rho} :=L_1(\mbb R, e^{\rho x}dx)$ for any $\rho>0$ (we restrict the considerations here
to the case $\sigma =0$ to avoid dealing with the singularity of the weight $x^\sigma$ at $x=0$).

As before, the solution to \eref{G5b} is given by
\begin{equation}
w^{+}(\xi, t) = (I+t\mc J^+)^m [\phi](\xi),\qquad \xi\in \mbb R,\qquad t\in [0,T],
\label{solw3}
\end{equation}
for  any $T<\infty$. Thus, we have
\begin{equation}
w^+(\xi, t) = (I+t\mc J^+)^m[v_0](\xi), \quad \xi>0
\label{sol15}
\end{equation}
and, using the definition of $\phi$, for $\xi<0$, we can write, see Appendix \ref{appB},
\begin{eqnarray}
w^+(\xi, t) &=&(I+t\mc J^+)^m[v_0](\xi) + (I+t\mc J)^m[\psi](\xi)\nonumber\\
& =:& F(\xi,t) + (I+t\mc J)^m[\psi](\xi),
\label{exsol1g}
\end{eqnarray}
where
\begin{equation}
\mc J[f](\xi):=\int_\xi^0 f(\eta)d\eta.
\label{mcJa}
\end{equation}
Thus, in view of \eref{cg1bcwb}, we get
 \begin{eqnarray}
 0&=&(I+t\mc J^+)^m[v_0](\xi)|_{t= -\beta^{-1}\xi} + (I+t\mc J)^m[\psi](\xi)|_{t= -\beta^{-1}\xi}\nonumber\\
& =& F(\xi, -\beta^{-1}\xi) + (I+t\mc J)^m[\psi](\xi)|_{t= -\beta^{-1}\xi}.
\label{exsol1f}
\end{eqnarray}
The problem is that the operators $(I+t\mc J^+)^m$ and $(I+t\mc J)^m$ are nonlocal and their evaluation
at a given $t$ for general $m$ is quite involved. Thus we shall proceed under the simplifying assumption that
$m \in \mbb N$.

Let us introduce the (probabilistic) Hermite polynomials, see \cite[Section 20.3]{abramowitz1964handbook},
\begin{equation}
He_m(\zeta) = m!\sum\limits_{i=0}^{\lfloor \frac{m}{2}\rfloor} \frac{(-1)^i}{i! (m-2i)!}\frac{\zeta^{m-2i}}{2^i}.
\label{Hek}
\end{equation}
Then, as shown in Appendix~\ref{appB}, the unique solution $\psi \in X^{+\rho}$ to \eref{exsol1f} reads
\begin{numparts}
\begin{eqnarray}
\psi(\xi) &=& (-1)^m\frac{d^m}{d\xi^m} \left(e^{-\frac{\xi^2}{2\beta}} y
\left(\frac{\xi}{\sqrt \beta}\right)\right)\nonumber\\
& =& \left.(-1)^m \beta^{\frac{m}{2}} \frac{d^m}{d\zeta^m}
\left (e^{-\frac{\zeta^2}{2}} y(\zeta)\right)\right |_{\zeta = \frac{\xi}{2}},
\label{psisol0}
\end{eqnarray}
 where $y$ is given by
 \begin{equation}
y(\zeta) = \int_0^\zeta \left(\sum\limits_{i=1}^m \frac{1}{(He_{m})'(\lambda_{m,i})}
e^{\lambda_{m,i}(\zeta-\sigma)}\right) e^{\frac{\sigma^2}{2}} g(\sigma) d\sigma,
\label{finsol0}
\end{equation}
\end{numparts}
$\lambda_{m,1}, \ldots,\lambda_{m,m}$ are simple real roots of $He_m$ and
$g(\zeta) = -\beta^{\frac{m}{2}}F\left(\zeta \beta^{\frac{1}{2}},-\beta^{-\frac{1}{2}}\zeta\right)$.

\subsection{Solutions}\label{sec4.4}
As mentioned earlier in Section~\ref{sec4.1}, the formulae for solutions in the constant growth/decay case
are based on the same representation \eref{solw11} as in Section \ref{sec3} and, consequently, their properties
can be established as in Appendix~\ref{App1}.

For the sake of completeness, we conclude this section by providing unified formulae for the solutions in
the decay case of \eref{eq3.1a}--\eref{eq3.1c}. Note that since we cover both positive and negative
$\alpha$, some restrictions below are superfluous.
The explicit solutions to \eref{eq3.1a}--\eref{eq3.1c} in the decay case are given by
\begin{eqnarray*}
u(x,t) & = 0,\quad x^\alpha < -k\alpha t,\quad x,t\in\mathbb{R}_+,\\
u(x,t) & = e^{-\frac{1}{2}ka\alpha t^2 - ax^\alpha t}
\Biggl[ \Bigl(1+\frac{k\alpha t}{x^\alpha}\Bigr)^{\frac{1-\alpha}{\alpha}}
u_0\Bigl(x \Biggl(1+\frac{k\alpha t}{x^\alpha}\Bigr)^{\frac{1}{\alpha}}\Biggr)\\
&+ a(\alpha+1) x^{\alpha-1} t \int_{(x^\alpha+k\alpha t)^{\frac{1}{\alpha}}}^{\infty}
{}_1F_1\Biggl(-\frac{1}{\alpha}; 2; at(x^\alpha+k\alpha t - y^\alpha) \Biggr)\\
&\qquad\qquad\qquad\qquad\qquad u_0(y)dy\Biggr],\quad x^\alpha \ge -k\alpha t,\quad x,t\in\mathbb{R}_+.
\end{eqnarray*}
Here, as $\mu = \theta = 0$, by default $\gamma = 1-\alpha$ and $\nu = \alpha-1$. In particular, for the
monodisperse initial data $u_0(x) = \delta(x-x_0)$, we have,
\begin{eqnarray*}
u(x,t) & = 0,\quad x_0^\alpha < k\alpha t,\quad x,t\in\mathbb{R}_+,\\
u(x,t) & = e^{-\frac{1}{2}ka\alpha t^2 - ax^\alpha t}\Biggl[
\delta\Biggl(x - x_0\Bigl(1-\frac{k\alpha t}{x_0^\alpha}\Bigr)^{\frac{1}{\alpha}}\Biggr)\\
&+ \chi_{[0,(x_0^\alpha-k\alpha t)^{\frac{1}{\alpha}}]}
a(\nu+2) x^{\alpha-1} t {}_1F_1\Biggl(-\frac{1}{\alpha}; 2; at(x^\alpha+k\alpha t - x_0^\alpha) \Biggr)\Biggr],\\
&x_0^\alpha \ge k\alpha t,\quad x,t\in\mathbb{R}_+.
\end{eqnarray*}

\begin{figure}[h!]
\begin{center}
\includegraphics{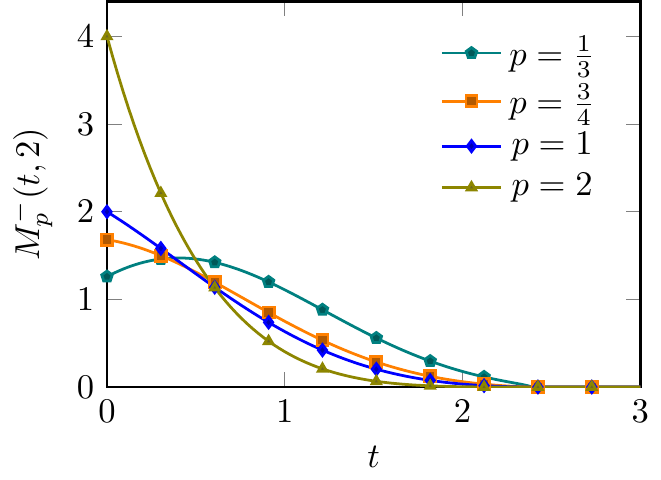} \hspace{2mm} \includegraphics{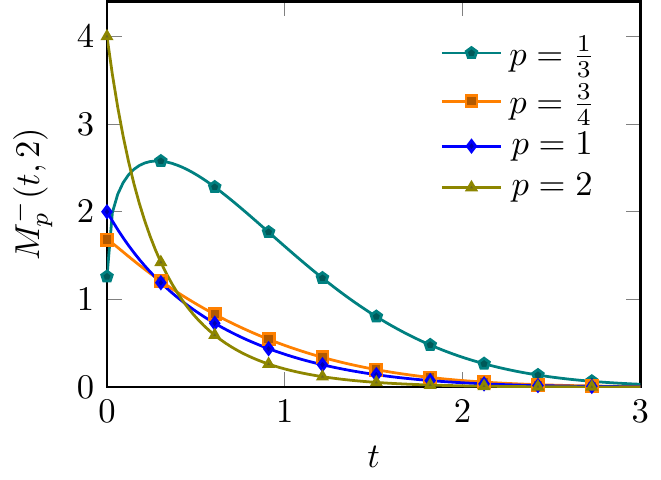}
\end{center}
\caption{Evolution of moments with $u^-_0(x) = \delta_2(x)$, $k=1$ and $a=1$:
left column $\alpha=\frac{2}{3}$; right column $\alpha=-\frac{3}{4}$.}\label{fig2}
\end{figure}

As in Section \ref{sec3}, we shall also provide formulae for the moments of solutions, focusing on the
monodisperse initial data only. First, we consider $\alpha>0$ and $p\ge 0$, in which case we have
\begin{eqnarray*}
M_p^-(t,x_0) & = 0, \quad x^\alpha < k\alpha t,\quad x_0,t\in\mathbb{R}_+,\\
M_p^-(t,x_0) & = e^{-\frac{1}{2}ka\alpha t^2} x_0^p\Bigl(1-\frac{k\alpha t}{x_0^\alpha}\Bigr)^{\frac{p}{\alpha}}
{}_1F_1\Biggl(\frac{p-1}{\alpha}; \frac{p+\alpha}{\alpha}; -at(x_0^\alpha-k\alpha t) \Biggr),\\
&x_0^\alpha \ge k\alpha t,\quad x_0,t\in\mathbb{R}_+.
\end{eqnarray*}
In particular,
\begin{eqnarray*}
M_1^-(t,x_0) & = 0, \quad x_0^\alpha < k\alpha t,\quad x_0,t\in\mathbb{R}_+,\\
M_1^-(t,x_0) & = e^{-\frac{1}{2}ka\alpha t^2} x_0\Bigl(1-\frac{k\alpha t}{x_0^\alpha}\Bigr)^{\frac{1}{\alpha}},
\quad x_0^\alpha \ge k\alpha t,\qquad x_0,t\in\mathbb{R}_+,
\end{eqnarray*}
and no shattering occurs.

On the other hand, when $-1<\alpha<0$ and $0<1+\alpha<p$, the moments are given by
\begin{eqnarray*}
M_p^-(t,x_0) & = \frac{\Gamma\bigl(\frac{\alpha-p+1}{\alpha}\bigr)}{\Gamma\bigl(\frac{p}{|\alpha|}\bigr)}
e^{\frac{1}{2}ka\alpha t^2-ax_0^\alpha t} x_0^p\Bigl(1-\frac{k\alpha t}{x_0^\alpha}\Bigr)^{\frac{p}{\alpha}}\\
&\phantom{xxxxx}\cdot\Psi\Biggl(\frac{\alpha+1}{\alpha}; \frac{\alpha+p}{\alpha}; at(x_0^\alpha-k\alpha t) \Biggr),
\qquad x_0,t\in\mathbb{R}_+,
\end{eqnarray*}
and, in particular, the formula
\begin{eqnarray*}
M_1^-(t,x_0) & = \frac{e^{-\frac{1}{2}ka\alpha t^2}}{\Gamma\bigl(\frac{1}{|\alpha|}\bigr)}
x_0\Bigl(1-\frac{k\alpha t}{x_0^\alpha}\Bigr)^{\frac{1}{\alpha}}
\Gamma\Bigl(\frac{1}{|\alpha|}; at(x_0^\alpha-k\alpha t) \Bigr),\quad
-1<\alpha<0,
\end{eqnarray*}
indicates that there exists spurious loss of mass not connected with the transport processes, i.e., for
$-1<\alpha<0$, we have shattering solutions.

The typical behavior of moments, with the monodisperse initial data $u^-_0(x)=\delta_{x_0}(x)$,
is shown in Fig~\ref{fig2}.

\section{Conclusion}
In this paper we considered continuous fragmentation equation with a transport term in the state space, describing either growth or decay of clusters. We have shown that in two special but important in application cases of linear or constant transport rates, the equations can be transformed to the same form as was obtained for the pure fragmentation equation in \cite{ZM1, ZiMc85, McZi87}. By interpreting the resulting equation as an ordinary differential equation with operator right hand side, we were able to write down solutions to all equations in a compact unified form as an algebraic function of the antiderivative operator acting on the initial condition. This operator form of the solution can be easily evaluated, giving explicit computable solutions in concrete cases. We have used these explicit solutions to gain physical insights into the described processes and validate the obtained solutions by evaluating the relevant moments of the solutions and comparing their properties with the theoretical predictions. In particular, we have shown that, indeed, shattering occurs in both growth and decay models if $\alpha<0$, in agreement with the results in \cite{6,Ed1} and \cite{BaAr}.

\appendix

\renewcommand{\thesection}{\Alph{section}}

\renewcommand{\numparts}{\addtocounter{equation}{1}%
     \setcounter{eqnval}{\value{equation}}%
     \setcounter{equation}{0}%
     \def\theequation{\ifnumbysec
     \Alph{section}.\arabic{eqnval}{\it\alph{equation}}%
     \else\arabic{eqnval}{\it\alph{equation}}\fi}}

\renewcommand{\endnumparts}{\def\theequation{\ifnumbysec
     \Alph{section}.\arabic{equation}\else
     \arabic{equation}\fi}%
     \setcounter{equation}{\value{eqnval}}}

\section{Justification of \eref{eq3.5a} and \eref{eq3.5b}}\label{App1}

\subsection{Extension of \eref{eq3.2}.}\label{secA.1}
It can be verified by direct calculations that $w^\pm$, with $\alpha<0$ are local, i.e.,
these solutions blow-up in a finite time in the sense of $X^{-\sigma}_{-\rho}$, $\sigma,\rho\in\mathbb{R}_+$.
For the forthcoming analysis, it is convenient to replace $w^\pm$ with
\begin{equation}\label{eq3.7}
f^\pm(\xi,\tau^\pm) := e^{-\tau^\pm\xi} w^{\pm}(\xi,\tau^\pm),\qquad
f^\pm_0(\xi):=w^\pm_0(\xi).
\end{equation}
It turns out that the latter functions are defined globally for $\tau^\pm\in\mathbb{R}_+$.
Furthermore, we have
\begin{lemma}\label{lm3.1}
Assume that $f_0^\pm\in X^{(\sg) \sigma}_0$ and either
$\alpha>0$ and $\sigma\ge 0$, $0\le \delta\le{m}$ or $\alpha<0$, $\sigma>{m}$
and $0\le \delta<\sigma-{m}$. Then
\begin{numparts}\label{eq3.8}
\begin{equation}\label{eq3.8a}
\nonumber
(\tau^{\pm})^{\delta}f^\pm \in C\bigl([0,T],X^{(\sg)\sigma+\delta}_0\bigr),
\end{equation}
for every finite value of $T>0$.
In addition, if $f_0^\pm\in X^{(\sg)\sigma+1}_0$, then
\begin{equation}\label{eq3.8b}
f^\pm \in C^1\bigl((0,T), X^{(\sg)\sigma+\delta}_0\bigr)
\end{equation}
\end{numparts}
and $f^\pm$ is the global classical solution to
\begin{numparts}\label{eq3.9}
\begin{eqnarray}
\label{eq3.9a}
& f^\pm_{\tau^\pm}(\xi,\tau^\pm) +\xi f^\pm(\xi,\tau^\pm) = m\mathcal{J}^{\sg}
\bigl[ f^\pm(\cdot,\tau^\pm) \bigr](\xi),
\qquad \xi,\tau^\pm\in\mathbb{R}_+,\\
\label{eq3.9b}
& f^\pm_{\tau^\pm}(\xi,0) = f_0^\pm(\xi),\qquad \xi\in\mathbb{R}_+,
\end{eqnarray}
\end{numparts}
in the sense of $X^{(\sg)\sigma}_0$.
\end{lemma}
\begin{proof}
(a) We let
\begin{eqnarray*}
&\mathcal{A}^\pm(\tau^\pm)[w_0^\pm](\xi) := e^{-\tau^\pm\xi} w_0^{\pm}(\xi),\\
&\mathcal{B}^\pm(\tau^\pm)[w_0^\pm](\xi) := m\tau^\pm  e^{-\tau^\pm\xi}\\
&\qquad\qquad\qquad
\mathcal{J}^{\sg}\Bigl[
{}_1F_1\Bigl(1-{(\sgn\alpha)m};2;-\tau^\pm(y-\xi) \Bigr) w^\pm_0\Bigr](\xi).
\end{eqnarray*}
Since $\xi^\delta e^{-\xi}\le c_\delta$, $\delta\in\mathbb{R}_+$, for some $c_\delta>0$,
uniformly in $\xi\in\mathbb{R}^+$, it follows that
\[
\|\mathcal{A}^\pm(\tau^\pm)\|_{X^{(\sg)\sigma}\to X^{(\sg)\sigma+\delta}}
\le c_\delta (\tau^\pm)^{-\delta},\qquad \sigma, \delta\ge 0,
\]
uniformly for $\tau^\pm\in\mathbb{R}_+$.  To estimate the norm of $\mathcal{B}^\pm$, we consider separately the cases $\alpha>0$ and $\alpha<0$.

(b) Assume initially that $\alpha>0$. Then, in view of \eref{eq3.3p}, we have
\begin{eqnarray*}
\|\mathcal{B}^\pm(\tau^\pm)\|_{X^{\sigma}\to X^{\sigma+\delta}}
&\le {m}(\tau^\pm)^{-\delta}\\
&\sup_{\xi\in\mathbb{R}_+} \xi^{-\sigma} e^{-\xi}
\int_0^\xi y^{\sigma+\delta} {}_1F_1\Bigl(1+{m};2;\xi-y \Bigr)dy.
\end{eqnarray*}
Next, by the asymptotic identity (see, e.g., \cite[Formulas 13.4.1 and 13.5.5, p.~508]{abramowitz1964handbook})
\begin{numparts}\label{eq3.10}
\begin{equation}\label{eq3.10a}
{}_1F_1\Bigl(a;b;-z \Bigr) =
\left\{
\begin{array}{ll}
\mathcal{O}(1),& z\to 0,\\
\mathcal{O}(z^{-a}),&z\to \infty,\quad \Rea z>0,
\end{array}
\right.
\end{equation}
and the formula
\begin{equation}\label{eq3.10b}
\int_0^x y^p {}_1F_1(a+1;2;x-y) dy
= \case{x^p}{a}
\left[{}_1F_1\Bigl(a;p+1; x  \Bigr)-1\right],
\end{equation}
\end{numparts}
which holds for all $a\ne0$ and $p>-1$, we have
\[
\|\mathcal{B}^\pm(\tau^\pm)\|_{X^{\sigma}\to X^{\sigma+\delta}}
\le c(\tau^\pm)^{-\delta},\qquad 0\le\delta\le{m},
\]
for some $c>0$. Hence, \eref{eq3.8a}, with $\alpha>0$, is settled.

(c) Let now $\alpha<0$. In this case,
\begin{eqnarray*}
\|\mathcal{B}^\pm(\tau^\pm)\|_{X^{-\sigma}\to X^{-\sigma+\delta}}
&\le {m}(\tau^\pm)^{-\delta}\\
&\sup_{\xi\in\mathbb{R}_+}
\xi^{\sigma} e^{-\xi}
\int_{\mathbb{R}_+} (\xi+y)^{\delta-\sigma} e^{-y}{}_1F_1\Bigl(1+{m};2;y \Bigr)dy.
\end{eqnarray*}
It is not difficult to verify that for $0<a<p$,
\begin{eqnarray}
\nonumber
&\int_{\mathbb{R}_+} (x+y)^{-p} e^{-y} {}_1F_1(a+1;2;y)dy
= \case{x^{-p}}{a\Gamma(p)} \int_{\mathbb{R}_+} e^{-t}t^{p-a-1}\Bigl[(x+t)^{a}-t^a\Bigr]dt\\
\nonumber
&\qquad\qquad\qquad
= \case{x^{-p}}{a\Gamma(p)}\Bigl[\Gamma(p-a)\Psi(-a;1-p;x)  - \Gamma(p)\Bigr]\\
\label{eq3.11}
&\qquad\qquad\qquad
\le
\left\{
\begin{array}{ll}
\case{\Gamma(p-a)}{a\Gamma(p)} x^{a-p},& 0<a\le 1,\\
\case{2^{a-1}\Gamma(p-a)}{a\Gamma(p)} x^{a-p} + \case{2^{a-1}-1}{a} x^{-p}, & a>1.
\end{array}
\right.
\end{eqnarray}
 Hence,
\[
\|\mathcal{B}^\pm(\tau^\pm)\|_{X^{-\sigma}\to X^{-\sigma+\delta}}
\le c_\delta(\tau^\pm)^{-\delta},\qquad 0\le\delta\le\sigma-{m}
\]
and the proof of \eref{eq3.8a} is complete.

(d) From \eref{eq3.10a}--\eref{eq3.10b}, \eref{eq3.11}, the inclusion $f_0^\pm\in X^{(\sgn\alpha)\sigma+1}_0$
and the standard identity (see \cite[Formula 13.4.8, p.~505]{abramowitz1964handbook})
\begin{eqnarray}\label{eq3.12}
\case{d}{dz}{}_1F_1(a;b;z) = \case{a}{b}{}_1F_1(a+1;b+1;z),
\end{eqnarray}
it follows (as in parts (b) and (c) above) that $f^\pm$, defined by the explicit formulas \eref{eq3.2} and \eref{eq3.7},
satisfies
\[
f^\pm_{\tau^\pm}, \xi f^\pm, \mathcal{J}^{\sg}[f^\pm]
\in C\bigl((0,T), X^{(\sg)\sigma+\delta}_0\bigr),
\]
for any finite value of $T>0$. Using this fact and the direct substitution of
$f^\pm$ into \eref{eq3.9a}--\eref{eq3.9b}, it is not difficult to verify that
\eref{eq3.9a} and \eref{eq3.9b} hold in $X_0^{(\sg)\sigma}$
globally for $\tau^\pm\in\mathbb{R}_+$.
\end{proof}

\subsection{Distributional and classical solutions.}\label{secA.2}
In connection with Lemma~\ref{lm3.1}, we note that $f^\pm$,
being integrable, satisfies \eref{eq3.9a}--\eref{eq3.9b} in the sense of Schwartz distributions. Moreover, from
\eref{eq3.3p} and \eref{eq3.12}, it follows that $f^\pm$, given by \eref{eq3.2} and \eref{eq3.7}, satisfies
\begin{eqnarray}
\nonumber
&&0=\int_{\mathbb{R}_+} f_0^\pm(\xi)v(\xi,0)
+\int_{\mathbb{R}_+}\left(\int_{\mathbb{R}_+} f^\pm(\xi,\tau^\pm)\right.\\
\label{eq3.13}
&&\phantom{xx}\left.\biggl[v_{\tau^\pm}(\xi,\tau^\pm) (\xi,\tau^\pm)
-\xi v(\xi,\tau^\pm)+ m\mathcal{J}^{-(\sgn\alpha)}[v(\cdot,\tau^\pm)](\xi)\biggr]
d\xi\phantom{\int_{\mathbb{R}_+} \!\!\!\!\!\!\!\!\!\!}\right)\! d\tau^\pm,
\end{eqnarray}
for any $v\in \mathcal{D}(\mathbb{R}_+^2)$ and $f^\pm_0\in \mathcal{D}'(\mathbb{R}_+)$. That is,
formulae \eref{eq3.2}, \eref{eq3.7} solve \eref{eq3.9a}--\eref{eq3.9b} in the sense of Schwartz distributions
for any distributional initial data in $\mathcal{D}'(\mathbb{R}_+)$.  In particular, for the monodisperse initial condition
$f_0^\pm(\cdot) = \delta_{\xi_0}(\cdot) = \delta(\cdot-\xi_0)$, supported at $\xi_0\in\mathbb{R}_+$,
we have
\begin{numparts}\label{eq3.14}
\begin{eqnarray}
\nonumber
f^\pm(\xi,\tau^\pm) &= e^{-\tau^\pm\xi}\Bigl[\delta_{\xi_0}(\xi) \\
\label{eq3.14a}
&+ \chi_{[0,\xi_0]}(\xi)  m\tau^\pm {}_1F_1\Bigl(1-{m};2; \tau^\pm(\xi-\xi_0) \Big)\Bigr],
\qquad \alpha>0,\\
\nonumber
f^\pm(\xi,\tau^\pm) &= e^{-\tau^\pm\xi}\Bigl[\delta_{\xi_0}(\xi) \\
\label{eq3.14b}
&+ \chi_{[\xi_0,\infty)}(\xi)  m\tau^\pm {}_1F_1\Bigl(1+{m};2;  \tau^\pm(\xi-\xi_0) \Big)\Bigr],
\qquad\alpha<0.
\end{eqnarray}
\end{numparts}
 As an immediate consequence
of the preliminary calculations, presented above, we have
\begin{corollary}\label{lm3.2}
For $u^\pm_0\in \mathcal{D}'(\mathbb{R}_+)$, the distributional solutions
to \eref{eq3.1a}--\eref{eq3.1b} are given explicitly by \eref{eq3.5a}.
In particular, for the monodisperse initial data $u^\pm_0 = \delta_{x_0}$, $x_0\in\mathbb{R}_+$,
formula \eref{eq3.5b} holds.

Further, for any finite value of $T>0$ and integrable input data $u^\pm_0\in X^p$,
with either $\alpha>0$, $p\ge \alpha-\nu-1$ and $0\le \delta\le \case{\nu+2}{\alpha}$
or $\alpha<0$, $p> 1+\alpha$ and $0\le \delta< \case{\alpha-p+1}{\alpha}$,
the solutions \eref{eq3.5a} satisfy
\begin{numparts}\label{eq3.16}
\begin{equation}\label{eq3.16a}
t^\delta u^\pm \in C\bigl([0,T],X^{p+\alpha\delta}_0\bigr),
\end{equation}
In addition, if the initial datum is regular, i.e., if $(u^\pm_0)_x\in X_0^{p+1}$ and $u^\pm_0\in X_0^{p+\alpha}$, then
\begin{eqnarray}
\label{eq3.16b}
&u^\pm_t, (ru^\pm)_x, au \in C\bigl((0,T),X^{p}_0\bigr),
\end{eqnarray}
\end{numparts}
and \eref{eq3.5a} satisfies \eref{eq3.1a}--\eref{eq3.1b} in the classical sense of $X_0^p$.
\end{corollary}
\begin{proof}
(a) By virtue of our definitions of $v^\pm$, $w^\pm$, $f^\pm$, $f_0^\pm$ and $\xi$, $z^\pm$, $\tau^\pm$,
the solution to \eref{eq3.1a}--\eref{eq3.1b} is formally given by
\begin{numparts}\label{eq3.17}
\begin{eqnarray}
\label{eq3.17a}
&u(x,t) = (1\pm\beta\tau^\pm(t))^{\frac{\nu}{\alpha}-\frac{\theta}{\beta}}
\xi^{\frac{\nu}{\alpha}}(x,t)f^\pm(\xi(x,t),\tau^\pm(t)),\\
\label{eq3.17b}
&\xi(x,t) = x^\alpha e^{\mp\beta t},\qquad
\tau^\pm(t) = \pm\case{1}{\beta}\bigl(e^{\pm\beta t}-1\bigr),\qquad x,t\in\mathbb{R}_+.
\end{eqnarray}
\end{numparts}
Since the coordinate transformation $(x,t)\mapsto(\xi,\tau^\pm)$, defined in \eref{eq3.17b}, is a
diffeomorphism from $\mathbb{R}_+^2$ to $\mathbb{R}_+ \times I^\pm$ and since
$f^\pm\in \mathcal{D}'(\mathbb{R}_+^2)$ satisfies \eref{eq3.9a}--\eref{eq3.9b} in the sense of distributions
for any $f_0^\pm\in \mathcal{D}'(\mathbb{R}_+)$, it follows (after changing variables in \eref{eq3.13})
that $u^\pm$, defined in \eref{eq3.5a}, satisfies
\begin{eqnarray*}
0=\int_{\mathbb{R}_+} u_0^\pm(x)v(x,0)dx
&+\int_{\mathbb{R}_+} dt\int_{\mathbb{R}_+} u^\pm(x,t)
\Biggl[ v_t(x,t)\pm r(x)v_x(x,t)\\
&- a(x)v(x,t) +a(x)\int_0^x b(y,x) v(y,t)dy
\biggr]dx,
\end{eqnarray*}
for any $u_0^\pm\in \mathcal{D}'(\mathbb{R}^+)$ and hence our first claim is settled.

(b) The right hand side of formula \eref{eq3.17a} defines one-to-one linear maps
$\mathcal{T}^\pm:f^\pm\mapsto u^\pm$. Elementary calculations shows that these maps satisfy
\begin{eqnarray*}
&\mathcal{T}^\pm\in L\bigl(C([0,T^\pm], X^\sigma_0), C([0,T], X^{p}_0)\bigr),\\
&\Bigr(\mathcal{T}^\pm\Bigl)^{-1}\in L\bigl(C([0,T^\pm], X^p_0), C([0,T], X^{\sigma}_0)\bigr),\\
&p =  \alpha(\sigma+1)-\nu-1,\qquad \sigma\in\mathbb{R},\qquad T = \pm\case{1}{\beta}(1\pm \beta T^\pm),
\end{eqnarray*}
for any finite $T^\pm\in I^\pm$. These inclusions, together with \eref{eq3.8a}--\eref{eq3.8b} and the identity
$\tau^\pm(t) = \mathcal{O}(t)$, $t\to0$, yield \eref{eq3.16a}.
In addition, if $(u^\pm_0)_x\in X_0^{p+1}$ and $u^\pm_0\in X_0^{p+\alpha}$,
direct calculations, using \eref{eq3.12}, \eref{eq3.17a}--\eref{eq3.17b} and \eref{eq3.8a}--\eref{eq3.8b},
show that $u^\pm$, defined in \eref{eq3.5a}, satisfy \eref{eq3.16b}.
Using this fact and the direct substitution, it is not difficult  to verify that $u^\pm$, defined in \eref{eq3.5a},
satisfy \eref{eq3.1a}--\eref{eq3.1b} in the classical sense of $X_0^p$. The proof is complete.
\end{proof}

\section{Proofs and additional formulae in the constant decay rate case}\label{appB}

\subsection{Proof of \eref{solw11}.}\label{secB.1}
\begin{proposition}\label{propB1}
The unique classical solution to \eref{G5aa}--\eref{G5ab} in $X^{\pm \sigma}_{\pm\rho},$ $\rho,\sigma>0$
is given by
\begin{equation}
w^{\pm}(\xi, t) = (I+t\mc J^\pm)^m [w_0](\xi),\qquad t\in [0,T],\qquad \xi\in \mathbb R_+,
\label{solw1}
\end{equation}
for any finite $T>0.$
\end{proposition}
\begin{proof} From Lemma \ref{lm2.1} we know that \eref{solw1} holds at least for $0< T<\rho$.
Further,
by directly solving the resolvent equation
\begin{equation}
\lambda f -t\mc  J^\pm f = g,
\label{reseq1}
\end{equation}
we find that the resolvent of $t\mc J^\pm$ is formally given by
\begin{equation}
(\lambda I-t\mc J^\pm)^{-1}[ g](\xi)
= \frac{t}{\lambda^2}\mc J^\pm [e^{\pm \frac{t}{\lambda}(\cdot - \xi)}g(\cdot)](\xi)
+\frac{1}{\lambda} g(\xi), \quad \xi\in \mathbb R_+.
\label{res1}
\end{equation}
Then, changing the order of integrals,
\begin{eqnarray*}
&&\int_0^\infty\xi^{\pm\sigma} e^{\pm \rho \xi}
\left|\mc J^\pm [e^{\pm \frac{t}{\lambda}(\cdot - \xi)}g(\cdot)](\xi)\right|d\xi\\
&&\le \left\{
\begin{array}{lcl}
\int_{0}^\infty \left(e^{\frac{t\Rea \lambda }{|\lambda|^2}\eta}| g(\eta)|
\int_{0}^\eta \xi^{\sigma}e^{\left(\rho-\frac{t\Rea \lambda }{|\lambda|^2}\right) \xi} d\xi \right)d\eta,
&\textrm{if}&\alpha>0,\\
\int_{0}^\infty \left(e^{-\frac{t\Rea \lambda }{|\lambda|^2}\eta}| g(\eta)|
\int_{\eta}^\infty \xi^{-\sigma}e^{-\left(\rho-\frac{t\Rea \lambda }{|\lambda|^2}\right) \xi} d\xi \right)d\eta,
&\textrm{if}&\alpha<0.
\end{array}
\right.
\end{eqnarray*}
Using monotonicity of $x^{\pm\sigma}$ on the respective intervals to factor it out from the inner integrals,
we see that if $\rho-\frac{t\Rea \lambda }{|\lambda|^2}>0$, then both integrals can be evaluated leading to
$$
\left\|\mc J^\pm [e^{\pm \frac{t}{\lambda}(\cdot - \xi)}g(\cdot)]\right\|_{X^{\pm \sigma}_{\pm \rho}}
\leq \frac{|\lambda|^2}{\rho|\lambda|^2-t\Rea\lambda}\|g\|_{X^{\pm \sigma}_{\pm \rho}}
$$
and hence the resolvent $(I-t\mc J^\pm)^{-1}$ exists in $X^{\pm \sigma}_{\pm \rho}$ as long as the above condition
is satisfied. Thus, the spectrum of $t\mc J^\pm$ is contained in
$\left\{ \lambda \in \mathbb C;\; \left(\Rea\lambda -\frac{t}{2\rho}\right)^2 + (\Ima \lambda)^2
\leq \frac{t^2}{4\rho^2}\right\}$,
which belongs to the closed right complex half-plane. Also, the spectrum of $0\mc J^\pm$ is 0.
Therefore, for any function $\mc F$ that is analytic in an open set containing the spectrum of
$t\mc J^\pm$ (for a fixed $t\geq 0$), we can evaluate $\Phi(t\mc J^\pm)$ by means of the Dunford integral
$$
\mc F(t\mc J^\pm) = \frac{1}{2\pi i}\int_\mathcal C \mc F(z)(I-t\mc J^\pm)^{-1} dz,
$$
where $\mathcal C$ is a curve surrounding the spectrum of $t\mc J^\pm$ in a positive direction.
We see that if we change $t\in [0,T]$, then the spectra of $t\mc J^\pm$ will continuously change from 0
to the disc centred at $\left(\frac{T}{2\rho},0\right)$ with radius $\frac{T}{2\rho},$
so each one will be contained in the latter. By the analyticity of the resolvent, we can define a smooth function
$[0,T]\ni t\mapsto \mc F(t\mc J^\pm)$ for any $0<T<\infty,$ provided the analyticity domain of $\mc F$ includes
the largest spectral disc of $t\mc J^\pm$. Since the functions $z \mapsto (1+tz)^m, m>0, t\in [0,T]$ are analytic in
$\mbb C$ with the cut along the negative ray  $\{z\in \mbb C;\; \Rea z> -1/T\}$ (with an obvious exception
for the integer values of $m$, which give polynomials), the solution (\ref{solw1}) can be extended to $[0,T]$
for any $T>0$.
\end{proof}

\paragraph{Remark.} The occurrence of the weight $\xi^{-\sigma}$ in the case $\alpha<0$
is natural if one has in mind transformation \eref{eq3.2p} for $z$.

\subsection{Proof of \eref{exsol1g}.}\label{secB.2}
To show that \eref{exsol1g} holds, we note that the formal resolvent is again given by \eref{res1},
with $\xi \in \mathbb R$ and,
as before, the spectrum of $t\mc J^+$ is contained in
$\left\{ \lambda \in \mathbb C;\; \left(\Rea\lambda -\frac{t}{2\rho}\right)^2 + (\Ima \lambda)^2
\leq \frac{t^2}{4\rho^2}\right\}.$ Hence, the solution to \eref{G5b} is given by the same formula
\eref{solw1}
$$
w^{+}(\xi, t) = (I+t\mc J^+)^m [\phi](\xi),\qquad t\in [0,T],
$$
for any $T<\infty$ but extended to $\xi\in \mbb R$.

For further use, we note that the solution to
 \begin{equation}
\lambda f(\xi) - t\int_\xi^0 f(\eta)d\eta =\lambda f(\xi) -t\mc J[f](\xi) = g(\xi), \qquad \xi\in \mathbb R_-,
\label{eqrestJ1}
\end{equation}
is given by
\begin{equation}
f(\xi) = \frac{t}{\lambda^2} e^{-\frac{t}{\lambda} \xi}
\int_\xi^0 e^{\frac{t}{\lambda}\eta} g(\eta)d\eta + \frac{1}{\lambda} g(\xi).
\label{restJ1}
\end{equation}
Using the Dunford integral representation for \eref{solw3}, we obtain
\begin{eqnarray}
w^+(\xi, t) &=& (I+t\mc J^+)^m[\phi](\xi)
= \frac{1}{2\pi i} \int_\mathcal C (1+z)^m (z I-t\mc J^+)^{-1}[\phi](\xi) dz\nonumber\\
&=&
\frac{1}{2\pi i} \int_\mathcal C (1+z)^m \left(\frac{t}{z^2} e^{-\frac{t}{z} \xi}
\int_\xi^\infty e^{\frac{t}{z}\eta} \phi(\eta)d\eta + \frac{1}{z} \phi(\xi)\right) dz.\label{exsol}
\end{eqnarray}
Now recall that for $\xi>0$ we have  $\phi(\xi)=v_0(\xi)$, with known $v_0$.  Thus, for such $\xi$, \eref{exsol}
provides a complete solution, given by (\ref{sol15}). Next, using the definition of $\phi$, for $\xi<0$ we can write
\begin{eqnarray}
w^+(\xi, t) &=&
\frac{1}{2\pi i} \int_\mathcal C (1+z)^m \left(\frac{t}{z^2} e^{-\frac{t}{z} \xi}
\int_0^\infty e^{\frac{t}{z}\eta} v_0(\eta)d\eta \right) dz\nonumber\\
&&+\frac{1}{2\pi i} \int_\mathcal C (1+z)^m \left(\frac{t}{z^2} e^{-\frac{t}{z} \xi}
\int_\xi^0e^{\frac{t}{z}\eta} \psi(\eta)d\eta + \frac{1}{z} \psi(\xi)\right) dz\nonumber\\
&=:& F(\xi,t) + (I+t\mc J)^m[\psi](\xi).
\label{exsol1}
\end{eqnarray}

\subsection{Proof of \eref{psisol0}--\eref{finsol0}.}\label{secB.3}
The proof of \eref{psisol0}--\eref{finsol0} is quite involved,  hence we split it into a series of lemmas.
\begin{lemma}\label{hermit}
The substitution $z = e^{-\frac{\zeta^2}{2}} y$ transforms the differential equation
\begin{equation}
\sum\limits_{k=0}^{m}\binom{m}{k} \zeta^{m-k} z^{(k)}(\zeta) =  g(\zeta), \qquad \zeta\in \mbb R,
\label{eqn2a}
\end{equation}
into to the constant coefficient equation
\begin{equation}
 He_m\left(\frac{d}{d\zeta}\right)[ y](\zeta) =  e^{\frac{\zeta^2}{2}} g(\zeta),
\label{eqn6}
\end{equation}
where $He_m$ is the probabilist's Hermite polynomial of order $m$, \cite[Section 20.3]{abramowitz1964handbook}.
\end{lemma}
\begin{proof}
Consider the substitution $z = h y$ for some known  differentiable function $h$.
Then, using the Leibnitz product formula and changing the order of summation,
\begin{eqnarray}
&&\sum\limits_{k=0}^{m}\binom{m}{k} \zeta^{m-k} (hy)^{(k)}(\zeta)
=\zeta^{m}\sum\limits_{r=0}^{m} y^{(r)}(\zeta)
\left(\sum\limits_{k=r}^m \binom{m}{k} \binom{k}{r} \zeta^{-k} h^{(k-r)} (\zeta) \right)\nonumber\\
&&=\zeta^{m}\sum\limits_{r=0}^{m} y^{(r)}(\zeta)
\left(\sum\limits_{l=0}^{m-r} \binom{m}{l+r} \binom{l+r}{r} \zeta^{-(l+r)} h^{(l)} (\zeta) \right)
\nonumber\\
&&=:\sum\limits_{r=0}^{m} a_{m,r}(\zeta) y^{(r)}(\zeta)=:L_m[y](\zeta).
\label{eqn3}
\end{eqnarray}
We shall prove that for $r\geq 1,$ we have
\begin{equation}
 a_{m+1, r}(\zeta) = \frac{m+1}{r} a_{m, r-1}.
 \label{acoeff}
 \end{equation}
 Indeed,
 \begin{equation}
 a_{m+1, r}(\zeta) = \zeta^{m+1}
\sum\limits_{l=0}^{m+1-r} \binom{m+1}{l+r} \binom{l+r}{r} \zeta^{-(l+r)} h^{(l)} (\zeta)
\label{am1}
\end{equation}
and
\begin{eqnarray*}
\binom{m+1}{l+r} \binom{l+r}{r}
&=& \frac{m+1}{r}\binom{m}{l+r-1}\binom{l+r-1}{r-1}.
\end{eqnarray*}
Thus, taking into account $\zeta^{m+1}\zeta^{-(l+r)} = \zeta^m \zeta^{-(l+(r-1))}$ in (\ref{am1}),
we obtain \eref{acoeff}.
Then, by iteration, \eref{acoeff} yields
\begin{equation}
a_{m,r}(\zeta) = \binom{m}{r} a_{m-r,0}(\zeta).
\label{amr}
\end{equation}
Hence, to specify all $a_{m,r}$, it suffices to know $a_{m,0}$ for any $m\in \mbb N_0$,
with $a_{0,0} = 1.$ In what follows,  we specify $h(\zeta) = e^{-\frac{\zeta^2}{2}}$
and note that
$$
h^{(k)}(\zeta) = h^{-1}(\zeta) (-1)^k He_k(\zeta),
$$
where $He_k$ is the Hermite polynomial, defined in \eref{Hek}.
Hence,
\begin{eqnarray*}
a_{m,0}(\zeta) 
&=&  h^{-1}(\zeta) \zeta^{m}\sum\limits_{l=0}^{m} \binom{m}{l}  (-1)^l
\left(\sum\limits_{i=0}^{\lfloor \frac{l}{2}\rfloor} \frac{(-1)^i}{i! (l-2i)!}\frac{\zeta^{-2i}}{2^i}\right).
\end{eqnarray*}
Considering the sum in the formula above, there are only even powers of $\zeta^{-1}$, running from $1$ to
$2\lfloor \frac{m}{2}\rfloor$.  The power $\zeta^{-2i}$ appears in the expansion only for $l\geq 2i$, with
the coefficient
\begin{eqnarray*}
b_{m,i} &:=&\sum\limits_{l=2i}^m \binom{m}{l}    \frac{(-1)^{l+i}l!}{i! (l-2i)!2^i}
= \frac{(-1)^i}{i!2^i}\sum\limits_{p=0}^{m-2i} (-1)^{p+2i}\frac{m!}{(m-p-2i)!p!} \\
&=& \frac{(-1)^i}{i!2^i}\prod\limits_{k=0}^{2i-1} (m-k)\sum\limits_{p=0}^{m-2i} (-1)^{p}\binom{m-2i}{p}.
\end{eqnarray*}
Thus, we have
$$
b_{m,i} = \left\{\begin{array}{lcl} 0 &\textrm{for}& 2i<m,\\
(-1)^\frac{m}{2} (m-1)!! &\textrm{for}& 2i=m.\end{array}
\right.
$$
Hence
\begin{equation}
a_{m,0}(\zeta) = e^{-\frac{\zeta^2}{2}} c_m := e^{-\frac{\zeta^2}{2}}
\left\{\begin{array}{lcl} 0 &\textrm{for}& m\;\textrm{odd},\\
(-1)^\frac{m}{2} (m-1)!! &\textrm{for}& m\;\textrm{even},\end{array}
\right.
\end{equation}
with the convention $(0-1)!! = 1$. Therefore, by \eref{amr},
\begin{eqnarray*}
&&L_m[y](\zeta)= e^{-\frac{\zeta^2}{2}}\sum\limits_{r=0}^{m} \binom{m}{r} c_{m-r} y^{(r)}(\zeta)
\\
&&=\left\{\begin{array}{lcl}
e^{-\frac{\zeta^2}{2}}\sum\limits_{i=0}^{k} \binom{2k}{2i}(-1)^{k-i}(2(k-i)-1)!! y^{(2i)}(\zeta)
&\!\!\textrm{for}& \!\!m = 2k,\\
e^{-\frac{\zeta^2}{2}}\sum\limits_{i=0}^{k} \binom{2k+1}{2i+1}(-1)^{k-i}(2(k-i)-1)!! y^{(2i+1)}(\zeta)
&\!\!\textrm{for}&\!\!m=2k+1.\end{array}
\right.
\end{eqnarray*}
Using the change of variable $l = k-i$ and the definition of the double factorial, the coefficients in both equations
can simplified to $(-1)^l\frac{m!}{(m-2l)!} \frac{1}{2^l l!}$ and, using \eref{Hek}, both differential operators can
be combined into
$$
L_m[y](\zeta) = e^{-\frac{\zeta^2}{2}}m!\sum\limits_{l=0}^{\lfloor \frac{m}{2}\rfloor}
\frac{(-1)^l}{l!(m-2l)!2^l } y^{(m-2l)}(\zeta) = e^{-\frac{\zeta^2}{2}}He_m\left(\frac{d}{d\zeta}\right)[ y](\zeta).
$$
Hence, \eref{eqn6} is proved.
\end{proof}

\begin{lemma}
The solution to \eref{eqn6} satisfying
\begin{equation}
y(0)=y'(0)=\cdots = y^{(m-1)} (0) =0,
\label{eqnic}
\end{equation}
 is given by
\begin{equation}
y(\zeta) = \int_0^\zeta \left(\sum\limits_{i=1}^m \frac{1}{(He_{m})'(\lambda_{m,i})}
e^{\lambda_{m,i}(\zeta-\sigma)}\right) e^{\frac{\sigma^2}{2}} g(\sigma) d\sigma,
\label{finsol}
\end{equation}
where  $\lambda_{m,1}, \ldots,\lambda_{m,m}$ are simple real roots of $He_m.$
\end{lemma}
\begin{proof}
The solution to (\ref{eqn6}) can be found by the variation of constants formula. Clearly, the characteristic
polynomial for (\ref{eqn6}) is $He_m(\lambda) = 0.$
By \cite[Section 22.16]{abramowitz1964handbook}, all zeroes of Hermite polynomials (as orthogonal polynomials)
are real and simple. Then the  functions $\theta_{m,i}(\zeta) = e^{\lambda_{m,i}\zeta}$ form a basis of the solution
space of the homogeneous equation \eref{eqn6}, hence we seek a particular solution to the inhomogeneous
equation \eref{eqn6} as
\begin{equation}
y(\zeta) = C_1(\zeta) e^{\lambda_{m,1}\zeta} +\cdots +  C_m(\zeta) e^{\lambda_{m,m}\zeta}.
\label{fory}
\end{equation}
In the general setting of the variation of constants method, $C_i'$s are given by
$$
C_i' = (-1)^{m-i}e^{\frac{\zeta^2}{2}} g(\zeta)\frac{W_i(\zeta)}{W(\zeta)},
$$
where $W$ is the Wronskian of $\{\theta_{m,1},\ldots,\theta_{m,m}\}$ and $W_i$ is the minor of the element $(m,i)$
of $W$. When the roots of the characteristic polynomial are simple, this formula can be made more explicit. Indeed,
\begin{eqnarray*}
W(\zeta) &=& 
\prod_{i=1}^{m}e^{\lambda_{m,i}\zeta}\left|\!\!\begin{array}{cccc} 1&1&\ldots&1\\
\lambda_{m,1}&\lambda_{m,2}&\ldots&\lambda_{m,m}\\
\vdots&\vdots&\vdots&\vdots\\
\lambda^{m-1}_{m,1}&\lambda^{m-1}_{m,2}&\ldots&\lambda^{m-1}_{m,m}\end{array}\!\!\right|=
V(\lambda_{m,1},\ldots,\lambda_{m,m})\prod_{i=1}^{m}e^{\lambda_{m,i}\zeta},
\end{eqnarray*}
where $V$ is the Vandermonde determinant, whose value is
$$
V(\lambda_{m,1},\ldots\lambda_{m,m}) = \prod_{1\leq i<j\leq m}(\lambda_{m,j}-\lambda_{m,i}).
$$
Now, for a given $r$,
$$
W_r = \prod_{i=1, i\neq r}^{m}e^{\lambda_{m,i}\zeta} V(\lambda_{m,1},\ldots,\lambda_{r-1},\lambda_{r+1},
\ldots,\lambda_{m,m}).
$$
To relate these two Vandermonde determinants, we write
\begin{eqnarray*}
&&V(\lambda_{m,1},\ldots\lambda_{m,m})\\
&&\phantom{xx} = \prod_{1<j\leq m}\!\!(\lambda_{m,j}-\lambda_{m,1})\cdot \ldots \cdot \prod_{r<j\leq m}
\!\!(\lambda_{m,j}-\lambda_{m,r})\cdot\ldots\cdot (\lambda_{m,m}-\lambda_{m,m-1})\\
&&=\phantom{xx} \prod_{1\leq i<r}(\lambda_{m,r}-\lambda_{m,i}) \cdot \prod_{r<j\leq m}
\!\!(\lambda_{m,j}-\lambda_{m,r}) \cdot\prod_{\stackrel{1\leq i<j\leq m}{ i,j\neq r}}\!\!(\lambda_{m,j}-\lambda_{m,i})\\
&&=\phantom{xx} (-1)^{m-r} \prod_{\stackrel{1\leq i\leq m}{i\neq r}}(\lambda_{m,r}-\lambda_{m,i})
\cdot V(\lambda_{m,1},\ldots,\lambda_{r-1},\lambda_{r+1},\ldots,\lambda_{m,m})\\
&&= (-1)^{m-r} (He_{m})'(\lambda_{m,r})V(\lambda_{m,1},\ldots,\lambda_{r-1},\lambda_{r+1},\ldots,\lambda_{m,m}).
\end{eqnarray*}
Hence, 
$$
C_i'(\zeta) = \frac{e^{\frac{\zeta^2}{2}} g_m(\zeta)}{e^{\lambda_{m,i}\zeta}(He_{m})'(\lambda_{m,i})}
$$
and \eref{finsol} follows from \eref{fory}.
\end{proof}

\begin{theorem}
Let $m\in \mathbb N$. Then the unique solution $\psi \in X^{+\rho}$ to \eref{exsol1} satisfying \eref{cg1bcwb}
is given by
\begin{eqnarray}
\psi(\xi) &=& (-1)^m\frac{d^m}{d\xi^m} \left(e^{-\frac{\xi^2}{2\beta}} y
\left(\frac{\xi}{\sqrt \beta}\right)\right)\nonumber\\
& =& \left.(-1)^m \beta^{\frac{m}{2}} \frac{d^m}{d\zeta^m}
\left (e^{-\frac{\zeta^2}{2}} y(\zeta)\right)\right |_{\zeta = \frac{\xi}{2}},
\label{psisol}
\end{eqnarray}
 where $y$ is given by \eref{finsol} with
$g(\zeta) = -\beta^{\frac{m}{2}}F\left(\zeta \beta^{\frac{1}{2}},-\beta^{-\frac{1}{2}}\zeta\right)$.
\end{theorem}
 \begin{proof} In this particular case, $F(\xi,t)$  is a known function, for $\xi<0$ given by
\begin{equation}
F(\xi,t) = \sum\limits_{r=1}^m\binom{m}{r}\frac{1}{(r-1)!} t^r\int_0^\infty (\eta-\xi)^{r-1}w_0(\eta) d\eta.
\label{Fm}
\end{equation}
Hence, we can re-write \eref{exsol1} as
\begin{equation}
(I+t\mc J)^m[\psi](\xi) = \sum\limits_{k=0}^{m}\binom{m}{k} t^k\mc J^k [\psi](\xi) = w^+(\xi,t) - F(\xi,t).
\label{eqn1}
\end{equation}
Now, we use \eref{cg1bcwb} and arrive at
\begin{equation}
\sum\limits_{k=0}^{m}\binom{m}{k} (-1)^k\beta^{-k}\xi^k\mc J^k [\psi](\xi) = - F(\xi,-\beta^{-1}\xi) =:G(\xi),
\label{eqn2}
\end{equation}
where $G$ is a known function. If we denote $Z(\xi) = (-1)^m\mc J^m [\psi](\xi),$ then
$Z^{(l)}(\xi) = (-1)^{m+l}\mc J^{m-l}[\psi](\xi)$ and \eref{eqn2} becomes
\begin{equation}
\sum\limits_{k=0}^{m}\binom{m}{k} \left(\frac{\xi}{\beta}\right)^{m-k} Z^{(k)}(\xi) = G(\xi).
\label{eqn23}
\end{equation}
Next, with $\zeta = \frac{\xi}{\sqrt \beta}$ and $z(\zeta) = Z(\xi)$,
\eref{eqn23} takes the form
\begin{equation}
\sum\limits_{k=0}^{m}\binom{m}{k} \zeta^{m-k} z^{(k)}(\zeta) = \beta^{\frac{m}{2}}G(\zeta \sqrt{\beta})
= : g(\zeta)
\label{eqn24}
\end{equation}
and hence, by Lemma \ref{hermit}, using $z(\zeta) = e^{-\frac{\zeta^2}{2}}y(\zeta)$, we transform \eref{eqn2a} into
\begin{equation}
He_m\left(\frac{d}{d\zeta}\right)[ y](\zeta)  = e^{\frac{\zeta^2}{2}}g(\zeta).
\label{herm}
\end{equation}
Further, we observe that, by definition, $Z^{k}(0) = 0$ for $k=0,\ldots, m-1$ and, since
$y(\zeta) = e^{-\frac{\zeta^2}{2}}z(\eta)$, the Leibniz formula shows that $y^{k}(0) = 0$, provided
$y^{(l)}(0) = 0$ for $l=0,\ldots,k-1,$ thus, by induction, the initial (or terminal) conditions for \eref{herm}
are given by \eref{eqnic}. Thus, the solution to \eref{eqn6} satisfying these conditions is given by \eref{finsol}
and we recover formula \eref{psisol} for $\psi$ by backward substitution.

To estimate $\psi$, we observe that, by \eref{Fm},  $e^{-\frac{\zeta^2}{2}} y(\zeta)$ is a linear combination of
terms of the form
$$
I_\lambda:= e^{\lambda \zeta - \frac{\zeta^2}{2}}
\int_0^\zeta p_{2m-1}(\sigma) e^{-\lambda \sigma + \frac{\sigma^2}{2}}d\sigma
= e^{-\frac{\upsilon^2}{2}} \int_{-\lambda}^{\upsilon} q_{2m-1}(\varsigma)
e^{\frac{\varsigma^2}{2}}d\varsigma,
$$
where $p_{2m-1}(\sigma)$ and $q_{2m-1}(\varsigma) = p_{2m-1}(\varsigma+\lambda)$ are polynomials of
degree $2m-1$ with coefficients depending on the moments of $w_0$ of order from $0$ to $m-1$.
Now, by the Leibniz rule, the $m$th derivative of $I_\lambda$ is a linear combination of $k$th derivatives of
$e^{-\frac{\upsilon^2}{2}}$ (which are the Hermite polynomials of degree $k$) and  $e^{-\frac{\upsilon^2}{2}}$
and the $(m-k)$th derivative of
$\int_{-\lambda}^{\upsilon} q_{2m-1}(\varsigma) e^{\frac{\varsigma^2}{2}}d\varsigma$, $0\leq k\leq m$, which,
apart from the case $k=m$, is the derivative of order $m-k-1$ of $q_{2m-1}(\varsigma) e^{\frac{\varsigma^2}{2}}.$
Using again the Leibniz rule, we see that the $r$th derivative of the latter is given by
$\bar q_{2m-1+r}(\varsigma)e^{\frac{\varsigma^2}{2}},$ where $\bar q_{2m-1+r}$ is a polynomial of degree
$2m-1+r$. Thus, the terms of $I_\lambda^{(m)}$ are products of polynomials of degrees $k$ and $3m-k-2$, and
hence polynomials of degree $3m-2$, except for $k=m,$ which  is given by
$$
He_m(\upsilon)e^{-\frac{\upsilon^2}{2}}\int_{-\lambda}^{\upsilon}
q_{2m-1}(\varsigma) e^{\frac{\varsigma^2}{2}}d\varsigma.
$$
By the l'H\^{o}spital rule,
$$
\lim\limits_{\upsilon\to -\infty} \frac{He_m(\upsilon)e^{-\frac{\upsilon^2}{2}}\int_{-\lambda}^{\upsilon}
q_{2m-1}(\varsigma) e^{\frac{\varsigma^2}{2}}d\varsigma}
{He_m(\upsilon)\bar {p}_{2m-2}(\upsilon)} = \lim\limits_{\upsilon\to -\infty}
\frac{\int_{-\lambda}^{\upsilon} q_{2m-1}(\varsigma) e^{\frac{\varsigma^2}{2}}d\varsigma}
{e^{\frac{\upsilon^2}{2}}\bar {p}_{2m-2}(\upsilon)} = \lim\limits_{\upsilon\to -\infty}
\frac{q_{2m-1}(\upsilon) }{\tilde {p}_{2m-1}(\upsilon)}
= l,
$$
for some finite $l$, where $\bar p_{2m-2}$ and $\tilde p_{2m-1}$ are polynomials of respective degrees.
Thus,
$$
\psi(\xi) = \mathcal O(\xi^{3m-2})\quad \textrm{as}\;\xi\to -\infty$$
and hence $\phi\in X^{+\rho}.$
\end{proof} 

\def\cprime{$'$} \def\cprime{$'$} \def\cprime{$'$} \def\cprime{$'$}
  \def\cprime{$'$} \def\cprime{$'$} \def\cprime{$'$} \def\cprime{$'$}
  \def\cprime{$'$} \def\cprime{$'$}
\providecommand{\newblock}{}

\end{document}